\documentclass[a4paper,10pt]{amsart}
\usepackage[foot]{amsaddr}
\usepackage{geometry}
\usepackage{graphicx}
\usepackage{lineno}
\usepackage{amsfonts} 
\usepackage{amsmath}
\usepackage{amssymb}
\usepackage{framed} 
\usepackage{multicol} 
\usepackage{tabularx}
\usepackage{tabulary}
\usepackage{siunitx}
\usepackage{caption}
\usepackage{float}
\usepackage{url} 
\usepackage{multicol}
\usepackage{breakurl} 
\usepackage[breaklinks]{hyperref} 
\usepackage{xcolor}
\usepackage{xpatch}
\usepackage{nomencl} 

\usepackage{amsmath}
\usepackage{amssymb}
\usepackage{mathrsfs}
\usepackage{color}
\usepackage{graphicx}
\usepackage[percent]{overpic}
\usepackage{url}
\usepackage{subfig}
\usepackage{enumerate}
\usepackage{amsthm}
\usepackage{pgf}
\usepackage{todonotes}
\usepackage{easyReview}

\def\f{{\bm f}}
\def\g{{\bm g}}
\def\n{{\bm n}}
\def\u{{\bm u}}

\def\0{\boldsymbol{0}}
\def\ss{\boldsymbol{\sigma}}

\def\dt{\partial_t}

\def\el {\nonumber }

\newcommand{\bm}[1]{\mbox{\boldmath{$#1$}}}
\def\ddiv{\nabla\cdot}

\makenomenclature
\makenomenclature 
\RequirePackage{ifthen}
\let\oldref\ref
\renewcommand{\ref}[1]{(\oldref{#1})}

 \graphicspath{
    {figures_ok/}
    {cyl/}
    {cavity/}
}

\usepackage{amssymb}

\DeclareMathAlphabet\mathbfcal{OMS}{cmsy}{b}{n}

\xpatchcmd{\thenomenclature}{%
  \section*{\nomname}
}{
}{\typeout{Success}}{\typeout{Failure}}

\usepackage{ifthen}
\renewcommand{\nomgroup}[1]{%
  \ifthenelse{\equal{#1}{A}}{\item[\textbf{Abbreviations}]}{%
    \ifthenelse{\equal{#1}{G}}{\item[\textbf{Symbols}]}{%
      \ifthenelse{\equal{#1}{C}}{\item[\textbf{Abbreviations}]}{%
        \ifthenelse{\equal{#1}{S}}{\item[\textbf{Subscripts}]}{%
          \ifthenelse{\equal{#1}{Z}}{\item[\textbf{Mathematical Symbols}]}{}
        }
      }
    }
  }
}

\ifpdf
\usepackage{pdfsync}
\fi

\usepackage{mathabx}

\begin{document}
\newcommand{\dif}{\mbox{d}}
\title[]{An efficient FV-based Virtual Boundary Method for the simulation of fluid-solid interaction}
\author{Michele Girfoglio\textsuperscript{1,*}}
\thanks{\textsuperscript{*}Corresponding Author.}
\address{\textsuperscript{1}SISSA, International School for Advanced Studies, Mathematics Area, mathLab Trieste, Italy.}
\email{mgirfogl@sissa.it}
\author{Giovanni Stabile\textsuperscript{1}}
\email{gstabile@sissa.it}
\author{Andrea Mola\textsuperscript{1}}
\email{amola@sissa.it}
\author{Gianluigi Rozza\textsuperscript{1}}
\email{grozza@sissa.it}


\date{}

\dedicatory{}


\begin{abstract} 
In this work, the Immersed Boundary Method (IBM) with feedback forcing introduced by \cite{Goldstein1993} and often referred in the literature as the Virtual Boundary Method (VBM), is addressed. 
The VBM has been extensively applied both within a Spectral and a Finite Difference (FD) framework. Here, we propose to combine the VBM with a computationally efficient Finite Volume (FV) method. We will show that for similar computational configurations, FV and FD methods provide significantly different results. Furthermore, we propose to modify the standard feedback forcing scheme, based on a Proportional-Integral (PI) controller, with the introduction of a derivative action, in order to obtain a Proportial-Integral-Derivative (PID) controller. The stability analysis for the Backward Differentiation Formula of order 1 (BDF1) time scheme is modified accordingly, and extended to the Backward Differentiation Formula of order 2 (BDF2) time scheme. We will show that, for the BDF2 time scheme, the derivative action allows to improve the stability characteristics of the system. Our approach is validated against numerical data available in the literature for a stationary/rigidly moving 2D circular cylinder in several configurations. Finally, a Fluid-Structure Interaction (FSI) benchmark, related to the frequency response of a cantilever beam coupled with a fluid, is presented: we numerically demonstrate that the introduction of the derivative action plays an important role in order to properly detect the fluid-structure interaction coupling.
\end{abstract}

\maketitle

\section{Introduction}\label{sec:intro}
In the scientific literature we may find works related with the Immersed Boundary Methods (IBM) \cite{Mittal2005, Kim2019} since early 1970's. Many researchers have progressively expanded their attention to these methods for their ability to simulate moving or deforming bodies characterized by complex surface geometries embedded in a fluid region. The main feature that makes IBM an useful and versatile technique is related to the fact that the Navier-Stokes equations are discretized over an orthogonal Cartesian grid and the effect of the no-slip boundary condition at the physical surface of the immersed body is obtained by introducing an additional body force term in the momentum equation.

The first IBM has been introduced by \cite{Peskin1972, Peskin1977, Peskin1981} in order to perform  numerical analyses for cardiovascular flows. 
The fluid flow is described by Eulerian variables defined on a fixed Cartesian mesh; on the other hand, the immersed boundary motion is described in a moving Lagrangian grid. In order to transfer variables between Eulerian and Lagrangian domains, a proper approximation of the Dirac delta function is used. 
Later, several authors have introduced modifications and extensions of the original method. Depending on whether the force is applied onto the continuous or discretized Navier-Stokes equations, these methods can be categorized into a continuous forcing approach and a discrete forcing approach, respectively \cite{Mittal2005}. 
The continuous forcing approach, in addition to the original Peskin method, includes also the Virtual Boundary Method (VBM) introduced by \cite{Goldstein1993}. Since the present work deals with the VBM, we focus on the literature review of this method.

\cite{Goldstein1993, Goldstein1995, Goldstein1998} developed and validated the VBM within a Spectral framework. Such works employ theory of control methodologies to impose that, in correspondance with the body location, the fluid velocity is coinciding with that of the solid. In particular, the method discussed employs a feedback forcing scheme, based on a Proportial-Integral (PI) controller, to enforce the no-slip condition at boundaries immersed in the fluid domain. Two free (negative) gains, one related to the proportional action and the other one related to the integral action, are used in the standard formulation of the virtual forcing. In these works, highly accurate spectral interpolation of the velocities from the grid points to the virtual boundary points was implemented whilst linear interpolation was used to distribute the effect of the forcing term to the nearby grid points. Examples including 2D flow past a cylinder, 3D turbulent channel and turbulent flow over riblets are presented. 
\cite{Saiki1996, Saiki1997} formulated an improved version of the method. The fluid velocities were here interpolated to the virtual boundary points by means of bilinear interpolation. On the other hand, the effect of the virtual boundary force was extrapolated back to the grid points by area-weighted averages. The method was coupled with high-order Finite Difference (FD) schemes so as to suppress the numerical oscillations caused by the forcing observed in the Chebyshev Spectral method used in \cite{Goldstein1993, Goldstein1995}. The method was used to simulate stationary, rotating, and oscillating cylinders in uniform flow as well as the transition process in a flat plate boundary layer. In these works \cite{Goldstein1993, Goldstein1995, Goldstein1998, Saiki1996, Saiki1997}, it was highlighted that the VBM is affected by very strict time-step restrictions since the amplitude of the feedback forcing needs to be large for accurately enforcing the boundary conditions, resulting in an extremely stiff system. \cite{Goldstein1993} performed the stability analysis when the forcing term is computed explicitly for a second-order Adams-Bashfort scheme. Very small Courant-Friedrichs-Lewy ($CFL$) numbers (of order $10^{-2}$) were reported in various cases, making the method unattractive. 
A partial and preliminary improvement of the stability limit was obtained by \cite{Fadlun2000} when the proportional action in the feedback forcing scheme is computed implicitly in time ($CFL \approx 10^{-1}$). Further improvements were provided by \cite{Lee2003} that extended the stability analysis performed by \cite{Goldstein1993} to the Backward Euler scheme (BDF1), third-order Adams-Bashfort schemes, and second-order and third-order Runge-Kutta schemes. Analytical predictions for the stability limits as functions of the forcing gains were provided and validated by using a Spectral method against three cases: 3D turbulent flow caused by a surface-mounted box, the flow around an impulsively starting cylinder, and the rough-wall turbulent boundary layer flow. It was shown that when the third-order Runge-Kutta scheme is adopted the VBM performs properly with $CFL \approx 1$. \cite{Shin2008} coupled the Peskin's regularized delta function approach with the VBM in order to relax the stability limits, i.e. to achieve large $CFL$ numbers, and to improve the transfer process between Eulerian and Lagrangian domains. The stability analysis process performed by \cite{Lee2003} was significantly simplified and precise stability limits related to the BDF1 scheme for 2-point, 3-point, 4-point and 6-point regularized delta functions were provided. The method was implemented in a Finite-Difference context and applied to the 2D flow past stationary and oscillating cylinders.  \cite{Fadlun2000, Iaccarino2003} showed that from a physical viewpoint the forcing scheme can be intepreted as a simple damped oscillator where the stiffness and the damping are related to the integral and proportional actions, respectively. Within this framework, \cite{Margnat2009} rigorously investigated the behaviour of the VBM as a second-order damped control system. The natural frequency and the damping coefficient are introduced as driving parameters of the method in the place of the usual gain coefficients. Reliable insights related to the role of each parameter as well as to the time-step optimisation were provided by considering simulations of flows involving sharp edges. It should be highlighted that the stability of the method depends not only on the values of free constans but also on the flow geometry \cite{Fadlun2000, Lee2003, Iaccarino2003, Margnat2009}. Thus, at the current state, no general rule to select the optimum values of the forcing gains is provided. We also report the contribution of \cite{Lee2006} that used the ``area-weighted" VBM \cite{Saiki1996} within a Finite-Volume framework: simulations related to 2D flow past a cylinder in several different configurations were provided. In conclusion, we mention contributions focused on special cases of the VBM obtained using only one of the two controller actions: see, e.g., \cite{Lai2000} (proportional controller) and \cite{Khadra2000} (integrating controller).

As discussed above, the VBM demonstrated to be efficient to simulate the presence of fixed or moving solid bodies immersed in a fluid domain. Moreover, the VBM has been successfull used within a Fluid-Structure Interaction (FSI) framework. See, e.g., \cite{Huang2007, Shin2010, Song2011, Qin2012, Uddin2015, Son2017} for examples. 
Recently \cite{Park2017} investigated the coupling between fluid-flexible body interactions and heat transfer by opening the door towards new multiphysics scenarios. In all these works, FD methods are used for the numerical discretization of both flow and structural governing equations.

As showed by the literature review, the VBM has been extensively investigated within a Finite Difference framework as well as a Spectral framework, although obiously other space approximations are possible. Here, we focus on Finite Volume (FV) methods, which have been used within discrete forcing approach IBMs \cite{Constant2017, Jasak2014}. To best of our knowledge, except for \cite{Lee2006}, the application of the VBM in a FV framework has been unexplored. Moreover, we observe that in \cite{Lee2006} some relevant computational parameters, such as time step size, $CFL$, number of Lagrangian points, error estimation and sensitivity analysis respect to the gain coefficients, are missing. In this manuscript, we intend to fill this gap by proposing an exhaustive analysis of the features of the VBM within a FV framework. We will show that for similar computational configurations (mesh refinement, time step size, time discretization scheme, gain coefficients, number of Lagrangian points) FV and FD methods provide significantly different results. Also, we propose to modify the classic feedback forcing scheme by introducing a derivative action in order to obtain a Proportial-Integral-Derivative (PID) controller as contemplated in \cite{Goldstein1993}. In order to highlight the role played by the derivative action, the stability analysis originally performed by \cite{Shin2008} for the BDF1 time scheme will be extended. Next, we will consider also the BDF2 time scheme that, to the best of our knowledge, has not been explored within the VBM framework. We will show that when the BDF1 scheme is used the derivative action deteriorates the stability characteristics of the system. On the contrary, when the BDF2 scheme is used, the derivative action improves the stability characteristics of the system and, in particolar, allows to obtain a stability region wider than the widest one related to the BDF1 scheme, that is obtained with the usual PI controller. Our approach is validated against numerical data available in the literature for a stationary/moving 2D circular cylinder in several configurations. Finally, we will present a FSI benchmark, related to the frequency response of a cantilever beam coupled with a surrounding fluid, where the introduction of the derivative action plays a crucial role in order to obtain reliable results.

All the computational results presented in this article have been performed with OpenFOAM\textsuperscript{\textregistered}\cite{Weller1998}, an open source finite volume C++ library widely used by commercial and academic organizations. See \cite{Constant2017, Jasak2014} for IBM techniques implemented in OpenFOAM\textsuperscript{\textregistered}. The FSI benchmark has been implemented by coupling OpenFOAM with an open source finite element C++ library, \emph{deal.ii} \cite{dealii}, within a partitioned approach framework. An important outcome of this work is that all the used for the preparation of this work are incorporated in open-source libraries\footnote{\url{https://mathlab.sissa.it/cse-software}} and therefore are readily available to the scientific community.

This work is organized as follows. In Sec. \ref{sec:problem_def}, we introduce the continuous formulation of the VBM. In Sec. \ref{sec:space_discrete}, we detail our strategy for space discretization, which combines the VBM with a Finite Volume method. The stability analysis is reported in Sec. \ref{sec:stability_analysis}, while numerical results are presented in Sec. \ref{sec:numerical_results}. Finally, conclusions and perspectives are drawn in Sec. \ref{sec:conclusion}.


\section{Problem definition}\label{sec:problem_def}

\subsection{The Navier-Stokes equations}
We consider the motion of an incompressible viscous fluid in a time-independent domain $\Omega$ 
over a time interval of interest $(t_0, T)$. The flow is described by the incompressible Navier-Stokes equations:
\begin{align}
\rho\,\dt \u + \rho\,\ddiv \left(\u \otimes \u\right) - \ddiv \ss & = \f\quad \mbox{ in }\Omega \times (t_0,T),\label{eq:ns-mom}\\
\ddiv \u & = 0\quad\, \mbox{ in }\Omega \times(t_0,T),\label{eq:ns-mass}
\end{align}
complemented by the boundary conditions
\begin{align}
\u & = \u_D\quad\ \mbox{ on } \partial\Omega_D \times(t_0,T),\label{eq:bc-ns-d}\\
\ss\cdot\n & = \g\quad \quad\mbox{ on } \partial\Omega_N \times(t_0,T),\label{eq:bc-ns-n}
\end{align}
and the initial data $\u = \u_0$ in $\Omega \times\{t_0\}$.
Here $\overline{\partial\Omega_D}\cup\overline{\partial\Omega_N}=\overline{\partial\Omega}$ and $\partial\Omega_D \cap\partial\Omega_N=\emptyset$. In addition 
$\rho$ is the fluid density, $\u$ is the fluid velocity, $\dt$ denotes the time derivative, $\ss$ is the Cauchy stress tensor, $\f$ is the momentum forcing applied to enforce the no-slip boundary condition along the immersed boundary, $\u_D,\g$ and $\u_0$ are 
given.
 Equation (\ref{eq:ns-mom}) represents the conservation of the linear momentum, while eq. (\ref{eq:ns-mass}) represents the conservation of the mass. For Newtonian fluids $\ss$ can be written as
\begin{equation}\label{eq:newtonian}
\ss (\u, p) = -p \mathbf{I} +\mu (\nabla\u + \nabla\u^T),
\end{equation}
where $p$ is the pressure and $\mu$ is the constant \emph{dynamic} viscosity.
Notice that by plugging \eqref{eq:newtonian} into eq.~\eqref{eq:ns-mom} this
can be rewritten as
\begin{align}
\rho\, \dt \u + \rho\,\ddiv \left(\u \otimes \u\right) - 2\mu\,\Delta\u  + \nabla p = \f\mbox{ in }\Omega \times (t_0,T).\label{eq:ns-lapls-1}
\end{align}

We define the Reynolds number as
\begin{equation}
\mbox{Re} = \frac{U_r L_r}{\nu}, \label{eq:re}
\end{equation}
where $\nu=\mu/\rho$ is the \emph{kinematic} viscosity of the fluid, and $U_r$ and $L_r$ are characteristic macroscopic velocity and length, respectively. 
\subsection{The feedback forcing scheme}
The interaction force between the fluid and the immersed boundary in the Lagrangian reference frame can be calculated by the feedback law \cite{Goldstein1993}:

\begin{equation}\label{eq:PI}
\bm F(s,t) = \alpha \int_0^t \left(\u_{ib}(s,\tau) - \u_{b}(s,\tau)\right) d \tau + \beta \left(\u_{ib}(s,t) - \u_{b}(s,t)\right),  
\end{equation}
where $s$ is the curvilinear abscissa, $\alpha$ and $\beta$ are large negative free constants, $\u_{ib}$ is the flow velocity obtained by the interpolation at the immersed boundary, and $\u_{b}$ is the velocity of the immersed boundary. The transfer of variables between Eulerian and Lagrangian domains is obtained by the Dirac delta function $\delta$ \cite{Peskin2002}:

\begin{equation}\label{eq:4}
\u_{ib}\left(s, t\right) = \int_\Omega \u\left(\bm{x},t\right) \delta \left(\bm{r}(s,t) - \bm{x}\right) d\bm{x}, \qquad 
\f\left(\bm{x}, t\right) = \int_\Gamma \bm{F}\left(s,t\right) \delta \left(\bm{x} - \bm{r}(s,t)\right) ds,
\end{equation}
where $\bm{r}(s,t)$ is the virtual boundary position and $\bm{x}$ are the locations of the nearby fluid grid points. The equation \eqref{eq:PI} provides a Proportional-Integral (PI) feedback control of the velocity near the immersed boundary \cite{Goldstein1993}. Also, it acts as a spring with damping, where the spring is represented by the integral control term (i.e., $\alpha$ is the stiffness), and the damping is represented by the proportional control term (i.e., $\beta$ is the damping coefficient) \cite{Fadlun2000, Iaccarino2003, Margnat2009}. Note that, to the best of our knowledge, all the works available in the literature focus on the feedback law \eqref{eq:PI} without a derivative action, although the idea to consider other control terms is contemplated in \cite{Goldstein1993}. In this work, we introduce this further contribution in order to obtain a Proportional-Integral-Derivative (PID) controller, 
\begin{equation}\label{eq:PID}
\bm F = \alpha \int_0^t \left(\u_{ib} - \u_{b}\right) d \tau + \beta \left(\u_{ib} - \u_{b}\right) + \gamma \left(\dot{\u}_{ib} - \dot{\u}_{b}\right).
\end{equation}
We observe that $\gamma$ can be interpreted as a coefficient that modifies the mass and inertia characteristics of the system. We will show that the derivative action strongly affects stability limits in the forcing parameters space (see Sec. \ref{sec:stability_analysis}) and plays an important part in FSI problems (see Sec. \ref{sec:FSI}). 



\section{Numerical discretization}\label{sec:space_discrete}

\subsection{Space discretization of Navier-Stokes equations: the Finite Volume approximation}

In this section we briefly discuss the space discretization of problems \eqref{eq:ns-lapls-1}-\eqref{eq:ns-mass}. 
We adopt the Finite Volume (FV) approximation that is derived directly from the integral form of the governing equations. We have chosen to implement the VBM within
the finite volume C++ library OpenFOAM\textsuperscript{\textregistered} \cite{Weller1998}. We partition the computational domain $\Omega$ into cells or control volumes $\Omega_i$, with $i = 1, \dots, N_{c}$, where $N_{c}$ is the total number of cells in the mesh. Let  \textbf{A}$_j$ be the surface vector of each face of the control volume. 

The integral form of eq.~\eqref{eq:ns-lapls-1} for each volume $\Omega_i$ is given by:

\begin{align}\label{eq:evolveFVtemp-1.1}
\rho \int_{\Omega_i} \dfrac{\partial \u}{\partial t} d\Omega + \rho\, \int_{\Omega_i} \ddiv \left(\u \otimes \u\right) d\Omega - 2\mu \int_{\Omega_i} \Delta\u d\Omega + \int_{\Omega_i}\nabla p d\Omega  = \int_{\Omega_i}{\bm f}d\Omega.
\end{align}
By applying the Gauss-divergence theorem, eq.~\eqref{eq:evolveFVtemp-1.1} becomes:

\begin{align}\label{eq:evolveFV-1.1}
\rho \int_{\Omega_i} \dfrac{\partial \u}{\partial t}d\Omega + \rho\, \int_{\partial \Omega_i} \left(\u \otimes \u\right) \cdot d\textbf{A} - 2\mu \int_{\partial \Omega_i} \nabla\u \cdot d\textbf{A} + \int_{\partial \Omega_i}p d\textbf{A}  = \int_{\Omega_i}{\bm f} d\Omega.
\end{align}
Each term in eq.~\eqref{eq:evolveFV-1.1} is approximated as follows: 

\begin{itemize}
\item[-] \textit{Gradient term}: 

\begin{align}\label{eq:grad}
\int_{\partial \Omega_i}p d\textbf{A} \approx \sum_j^{} p_j \textbf{A}_j, 
\end{align}
where $p_j$ is the value of the pressure relative to centroid of the $j^{\text{th}}$ face. 
The face center pressure values $p_j$ are obtained from the cell center values by means of a linear interpolation scheme. 


\item[-] \textit{Convective term}: 
\begin{align}\label{eq:conv}
\int_{\partial \Omega_i} \left(\u \otimes \u\right) \cdot d\textbf{A} \approx \sum_j^{} \left(\u_j \otimes \u_j\right) \cdot \textbf{A}_j = \sum_j^{} \varphi_j \u_j, \quad \varphi_j = \u_j \cdot \textbf{A}_j,
\end{align} 
where $\u_j$ is the fluid velocity relative to the centroid of each control volume face. In \eqref{eq:conv}, $\varphi_j$ is the convective flux associated to $\u$ through face $j$ of the control volume. The convective flux at the cell faces is obtained by a linear interpolation of the values from the adjacent cells. Also $\u$ needs to be approximated at cell face $j$ in order to get the face value $\u_j$. Different interpolation methods can be applied: central, upwind, second order upwind and blended differencing schemes \cite{jasakphd}. In this work, we make use of a second order upwind scheme.
\item[-] \textit{Diffusion term}: 
\begin{align}
\int_{\partial \Omega_i} \nabla\u \cdot d\textbf{A} \approx \sum_j^{} (\nabla\u)_j \cdot \textbf{A}_j, \el
\end{align} 
where $(\nabla\u)_j$ is the gradient of $\u$ at face $j$. 
We are going to briefly explain how $(\nabla\u)_j$ is approximated with
second order accuracy on structured, orthogonal meshes, that are used in this work. Let $P$ and $Q$ be two neighboring control volumes.
The term $(\nabla\u)_j$ is evaluated by subtracting
 the value of velocity at the cell centroid on the $P$-side of the face, denoted with $\u_P$,
 from the value of velocity at the centroid on the $Q$-side, denoted with $\u_Q$,
 and dividing by the magnitude of the distance vector $\textbf{d}_j$ connecting the two cell centroids:
\begin{align}
(\nabla\u)_j \cdot \textbf{A}_j = \dfrac{\u_Q - \u_P}{|\textbf{d}_j|} |\textbf{A}_j|. \el
\end{align} 
\end{itemize}

A partitioned approach has been used to deal with the pressure-velocity coupling. In particular a Poisson equation for pressure has been used. This is obtained by taking the divergence of the momentum equation \eqref{eq:ns-lapls-1} and exploiting the divergence free constraint \eqref{eq:ns-mass}:
\begin{equation}\label{eq:Poisson}
\Delta p = -\nabla \left(\u \otimes \u\right).
\end{equation}

The segregated algorithms available in OpenFOAM\textsuperscript{\textregistered} are SIMPLE \cite{SIMPLE} for steady-state problems, and PISO \cite{PISO} and PIMPLE \cite{PIMPLE} for transient problems. For this work, we choose the PISO algorithm. 

\subsection{Time discretization}
To discretize in time the equation \eqref{eq:evolveFV-1.1}, let $\Delta t \in \mathbb{R}$, $t^n = t_0 + n \Delta t$, with $n = 0, ..., N_T$ and $T = t_0 + N_T \Delta t$. Moreover, we denote by $\u^n$ the approximation of the flow velocity at the time $t^n$. We adopt Backward Differentiation Formula of order 1 (BDF1) and Backward Differentiation Formula of order 2 (BDF2), see e.g. \cite{quarteroni2007numerical}. Given $\u^n$, for $n \geq 0$, we have, respectively,

\begin{equation}\label{eq:BDF1_disc}
\partial_t \u \approx \dfrac{\u^{n+1} - \u^{n}}{\Delta t},
\end{equation}

\begin{equation}\label{eq:BDF2_disc}
\partial_t \u \approx \dfrac{3\u^{n+1} - 4\u^{n} + \u^{n-1}}{2\Delta t}.
\end{equation}

\subsection{The interaction force} 
We represent the virtual boundary as a discrete set of $N_L$ Lagrangian points, $k$, with $k = 0, \dots, N_L$. The fluid-solid interaction force term is integrated over time using an explicit scheme: at the time $t^n$, we have

\begin{equation}\label{eq:forc1}
\bm F_k^{n} = \alpha \sum_{l = 1}^{n} ((\u_{ib})_k^{l} - (\u_{b})_k^{l})\Delta t + \beta ((\u_{ib})_k^{n} - (\u_{b})_k^{n}) + \gamma ((\dot{\u}_{ib})_k^{n} - (\dot{\u}_{b})_k^{n})
\end{equation}
where $(\u_{b})_k^{l}$ and $(\u_{ib})_k^{l}$, $(\dot{\u}_{b})_k^{l}$ and $(\dot{\u}_{ib})_k^{l}$, are computed as follows

\begin{equation}\label{eq:forc2}
(\u_{b})_k^{l} = \dfrac{\bm{r}_k^l - \bm{r}_k^{l-1}}{\Delta t}, \qquad (\u_{ib})_k^{l} = \sum_{\bm{x}_m \in \tilde{g}}  \u_m^l \tilde{\delta}(\bm{r}_k^l - \bm{x}_m)\Omega_m,
\end{equation}

\begin{equation}\label{eq:forc2_bis}
(\dot{\u}_{b})_k^{l} = \dfrac{(\u_b)_k^l - (\u_b)_k^{l-1}}{\Delta t}, \qquad (\dot{\u}_{ib})_k^{l} = \dfrac{(\dot{\u}_{ib})_k^l - (\dot{\u}_{ib})_k^{l-1}}{\Delta t},
\end{equation}
$\tilde{g}$ is the support of the smoothed delta function $\tilde{\delta}$ defined as

\begin{equation}\label{eq:forc3}
\tilde{\delta}(\bm{x}) = \dfrac{1}{\Omega_m}\prod_{\omega=1}^2 \varphi \left(\dfrac{x_{\omega}}{h}\right),
\end{equation}
$h$ the uniform mesh size next to the immersed boundary and $\Omega_m$ is given by
\begin{equation}\label{eq:vol2D}
\Omega_m = h^2.
\end{equation}

In this work, we use the 4-point regularized delta function \cite{Peskin2002, Shin2008}: 

\begin{equation}
\varphi(r) = 
\begin{cases}
\dfrac{1}{8}\left(3 - 2|r| + \sqrt{1+4|r|-4r^2}\right), & \text{if $0 \leq |r| \leq 1$}\\ \\
\dfrac{1}{8}\left(5 - 2|r| - \sqrt{-7+12|r|-4r^2}\right), & \text{if $1 \leq |r| \leq 2$}\\ \\
0, & \text{otherwise}.
\end{cases}
\end{equation}
Finally, the interaction force in the Eulerian reference frame is given by

\begin{equation}\label{eq:forc5}
\bm f_m^{n} = \sum_{k=1}^{N_L}  \bm F_k^{n} \tilde{\delta}(\bm{x}_m - \bm{r}_k^n) \Delta V,
\end{equation}
where $\Delta V$ is defined as \cite{Shin2008, Uhlmann2005}

\begin{equation}\label{eq:forc6}
\Delta V = h \cdot \Delta s,
\end{equation}
and $\Delta s$ is the uniform distance between two Lagrangian points.






\section{Stability analysis}\label{sec:stability_analysis}
In this section, the stability analysis of the VBM is presented. Concerning the standard formulation of the feedback forcing based on the PI controller, i.e. with $\gamma = 0$, precise stability boundaries in the forcing parameters space related to several time schemes were discussed in previous works. In particular, BDF1 \cite{Lee2003, Shin2008}, second-order and third-order Runge-Kutta schemes \cite{Lee2003}, second-order and third-order Adams-Bashfort schemes \cite{Goldstein1993, Lee2003}. Here, we are going to analyze the influence of the derivative controller term on the stability features of the BDF1 scheme. Furthermore, we will extend the analysis to the BDF2 time scheme that, to the best of our knowledge, has not been explored within the VBM framework. 

Following the approach performed in \cite{Shin2008}, the maximum magnitude of the Eulerian forcing can be expressed as:
\begin{equation}\label{eq:max_forc}
f_{max}(t) = C_{max} F(t) = C_{max} \left(\alpha \int_0^t \left(\u - \u_{b}\right) d \tau + \beta \left(\u - \u_{b}\right) + \gamma \left(\dot{\u} - \dot{\u}_{b}\right)\right),
\end{equation}
where $C_{max}$ is a coefficient depending on the type of regularized delta function as well as the dimensionality of the problem \cite{Shin2008}. 
For the sake of simplicity, we set $\u_{b} = \dot{\u}_{b} = 0$. Since $\alpha$ and $\beta$ makes the forcing term much larger than the other terms of the equation \eqref{eq:ns-lapls-1}, we can limit to investigate the stability features of the following equation \cite{Shin2008}:
\begin{equation}\label{eq:ODE1}
\dfrac{d \u}{dt} = C_{max} \left(\alpha \int_0^t \u d \tau + \beta \u + \gamma \dot{\u}\right).
\end{equation}






\subsection{BDF1 time scheme}

When the BDF1 scheme is adopted for the temporal discretization, we have

\begin{equation}\label{eq:BDF1_classic}
\u^{n+1} - \u^{n}= \alpha^\star \sum_{i=0}^{n} \u^{i} + \beta^\star \u^n + \gamma^\star \left(\u^n - \u^{n-1}\right),
\end{equation}
where $\alpha^\star = C_{max}\alpha \Delta t^2$, $\beta^\star = C_{max}\beta \Delta t$ and $\gamma^\star = C_{max}\gamma$. 

In order to obtain the recurrence formula for the stability analysis, the equation at the previous time step, $n-1$, is subtracted from the equation at the present time step, $n$, resulting in:
\begin{equation}\label{eq:BDF1_classic_2}
\u^{n+1} - 2\u^{n} + \u^{n-1}  = \alpha^\star \u^{n} + \beta^\star \left( \u^{n} - \u^{n-1}\right) + \gamma^\star  \left(\u^{n} - 2\u^{n-1} + \u^{n-2}\right).
\end{equation}

The stability features are related to $r = \u^{i+1}/\u^{i}$, $i = 0,1,2,...,n$. By recasting the equation \eqref{eq:BDF1_classic_2} in terms of $r$, we obtain
\begin{equation}\label{eq:BDF1_stab}
r^3 - (2 + \alpha^\star + \beta^\star + \gamma^\star)r^2 + (1 + \beta^\star + 2 \gamma^\star) r - \gamma^\star = 0.
\end{equation}

Notice that, for $\gamma^\star = 0$, we obtain the characteristic equation reported in \cite{Shin2008}:
\begin{equation}\label{eq:BDF1_stab2}
r^2- (2 + \alpha^\star + \beta^\star + \gamma^\star)r + 1 + \beta^\star = 0.
\end{equation}

The stability region, $|r|\leq1$, can be found by using the Jury criterion. For $\gamma^\star = 0$, we have \cite{Shin2008}:
\begin{equation}\label{eq:stab1}
-\alpha^\star - 2\beta^\star \leq 4.
\end{equation}

Now, we are going to investigate the influence of $\gamma^\star$ on the stability region. Since the stability region is too complicated to be expressed in closed form by leaving free all the parameters, $\alpha^\star$, $\beta^\star$ and $\gamma^\star$, we fix the value of $\gamma^\star$ and limit to investigate some specific cases in order to obtain a qualitative trend. Based on the equation \eqref{eq:BDF1_stab}, necessary condition for the stability is that $\gamma^\star \in \left[-1, 1\right]$. We are only interested in investigating negative values. 
We set $-\gamma^\star = 0.25$, $-\gamma^\star = -0.5$, $-\gamma^\star = 0.75$, and $-\gamma^\star = 1$, by obtaining, respectively, 
\begin{equation}\label{eq:stab2}
\begin{cases}
-\alpha^\star - 2\beta^\star \leq 3, \\
-\alpha^\star - 2\beta^\star \leq 2, \\
-\alpha^\star - 2\beta^\star \leq 1, \\
-\alpha^\star = -\beta^\star = 0.
\end{cases}
\end{equation}

Stability regimes \eqref{eq:stab1}-\eqref{eq:stab2} in 2D flow, obtained by set $C_{max} = \dfrac{1}{2}$ \cite{Shin2008}, 
\begin{equation}\label{eq:stab2_without_star}
\begin{cases}
-\alpha \Delta t^2 - 2\beta \Delta t \leq 8 & \text{for $-\gamma = 0$}, \\
-\alpha \Delta t^2- 2\beta \Delta t \leq 6 & \text{for $-\gamma = 0.5$}, \\
-\alpha \Delta t^2- 2\beta \Delta t \leq 4 & \text{for $-\gamma = 1$}, \\
-\alpha \Delta t^2 - 2\beta \Delta t \leq 2 & \text{for $-\gamma = 1.5$}, \\
-\alpha \Delta t^2 = -\beta \Delta t = 0 & \text{for $-\gamma = 2$},
\end{cases}
\end{equation}

are displayed in Figure \ref{fig:BDF_stability_1} a). The flow is stable in the region below the line and unstable above the line. We observe that the widest stability region is obtained for $\gamma = 0$, its size reduces, i.e. the intercept of the line moves downwards, to increasing of $-\gamma$ and finally degenerates into the axis origin for $-\gamma = 2$. On the other hand, the shape of the stability region is not affected by $\gamma$. In order to validate such theoretical predictions, simulations of transverse oscillations of a circular cylinder in a free-stream (see Sec \ref{sec:osc_y_cyl}) are performed for $\gamma$ = 0 (Figure \ref{fig:BDF_stability_2} a)) and $-\gamma$ = 1 (Figure \ref{fig:BDF_stability_3} a)). Numerical stable and unstable cases are denoted by circles and crosses, respectively. We observe that the numerical stability regions are wider than analytical ones because the analytical predictions are obtained based on the maximum value of the forcing (eq. \ref{eq:max_forc}). This result is in agreement with \cite{Shin2008}. For $\gamma$ = 0, we observe that the simulation is stable until to the line $-\alpha\Delta t^2 - 2\beta\Delta t = 11$. On the contrary, in \cite{Shin2008}, the simulation is stable until to the curve $-\alpha\Delta t^2 - 2\beta\Delta t = 10$. However, the marginal stability line is expected to change depending on the kind of flow.

In conclusion, we learned that when the BDF1 time scheme is used, the PI controller is to be preferred to the PID controller in order to obtain larger stability regions. Therefore, the derivative action deteriorates the stability features of the system. 

\subsection{BDF2 time scheme}

For the BDF2 time scheme, the approximation of the equation \eqref{eq:ODE1} yields
\begin{equation}
3\u^{n+1} - 4\u^{n} + \u^{n-1} = 2\left(\alpha^\star \sum_{i=0}^{n} \u^{i} + \beta^\star \u^n + \gamma^\star \left(\u^n - \u^{n-1}\right)\right).
\end{equation}
Notice that for $-\gamma^\star = 0.5$ the solution degenerates into that one obtained for the BDF1 time scheme with $\gamma^\star = 0$ unless a scale factor. 
The characteristic equation determining $r$ for this scheme is given by
\begin{equation}\label{eq:stab_BDF2}
3r^3 - (7 + 2 \alpha^\star + 2 \beta^\star + 2\gamma^\star)r^2 + (5 + 2 \beta^\star + 4 \gamma^\star) r - (1+2\gamma^\star) = 0.
\end{equation}

For $\gamma^\star = 0$, the stability region $|r|\leq1$ is

\begin{equation}\label{eq:stab1BDF2}
\begin{cases}
- \alpha^\star - 2 \beta^\star \leq 8, \\
 \alpha^\star - 2 \beta^\star \geq 0.\\
\end{cases}
\end{equation}

As we did for the BDF1 scheme, we investigate the influence of $\gamma^\star$ on the stability region. Based on the equation \eqref{eq:stab_BDF2}, necessary condition for the stability is $\gamma^\star \in \left[-2, 1\right]$. We consider only negative values. We set $\gamma^\star = -0.4$,  $\gamma^\star = -0.5$, $\gamma^\star = -0.6$, $\gamma^\star = -1$, $\gamma^\star = -1.5$, and $\gamma^\star = -2$, by obtaining, respectively,

\begin{equation}\label{eq:stab2BDF2}
\begin{cases}
\begin{cases}
- \alpha^\star - 2 \beta^\star \leq 6.4, \\    
 \alpha^\star - 14 \beta^\star \geq 0,\\
\end{cases}\\
-\alpha^\star - 2\beta^\star \leq 6, \\
-\alpha^\star - 2\beta^\star \leq 5.6, \\
-\alpha^\star - 2\beta^\star \leq 4, \\
-\alpha^\star - 2\beta^\star \leq 2,\\
-\alpha^\star = -\beta^\star = 0. \\
\end{cases}
\end{equation}

\begin{figure}
\centering
 \begin{overpic}[width=0.45\textwidth]{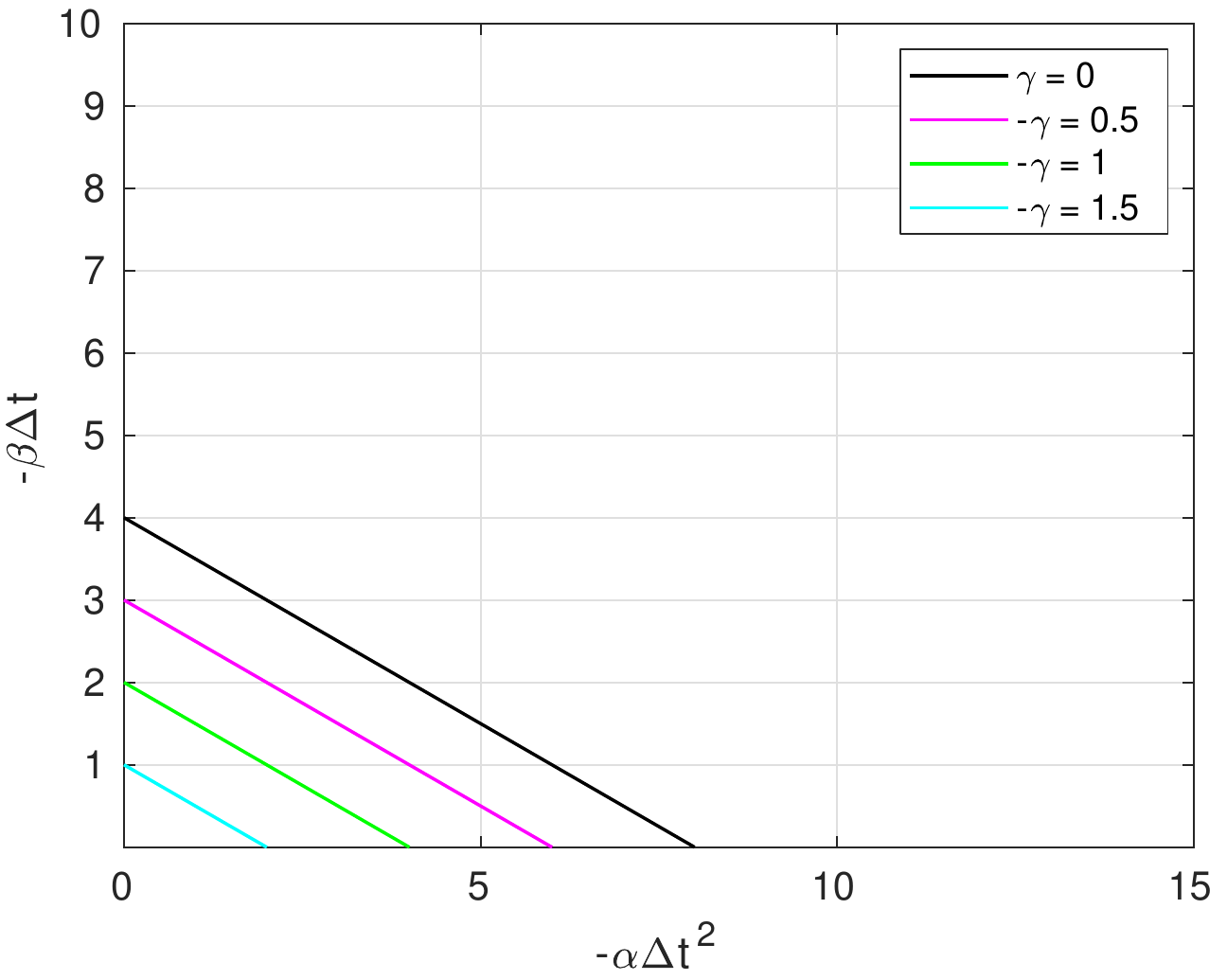}
       \put(50,82){\small{a)}}
      \end{overpic}
 \begin{overpic}[width=0.45\textwidth]{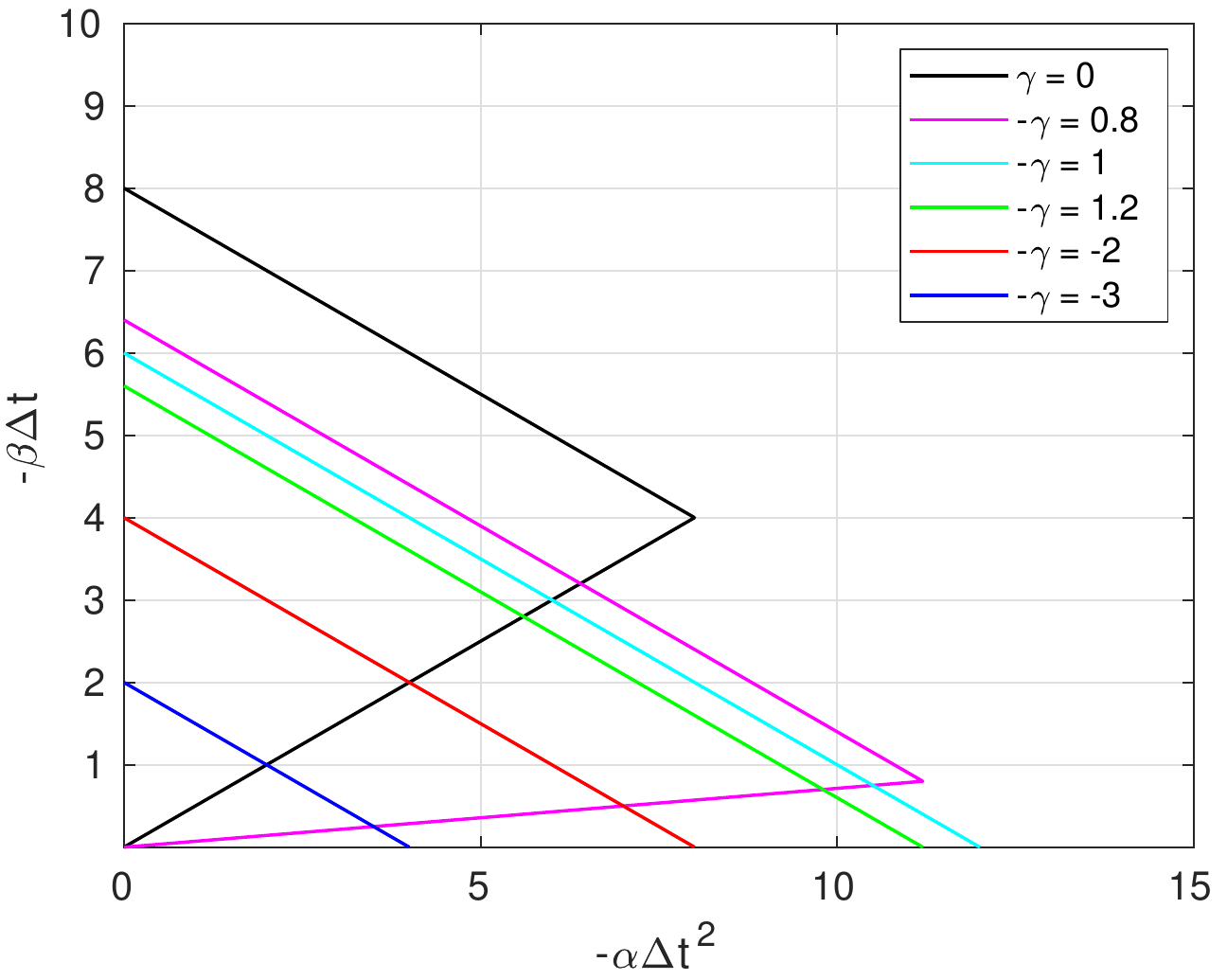}
       \put(50,82){\small{b)}}
   \end{overpic}
\caption{Analytical stability regions in 2D flow: a) BDF1 time scheme, b) BDF2 time scheme.}
\label{fig:BDF_stability_1}
\end{figure}

\begin{figure}
\centering
 \begin{overpic}[width=0.45\textwidth]{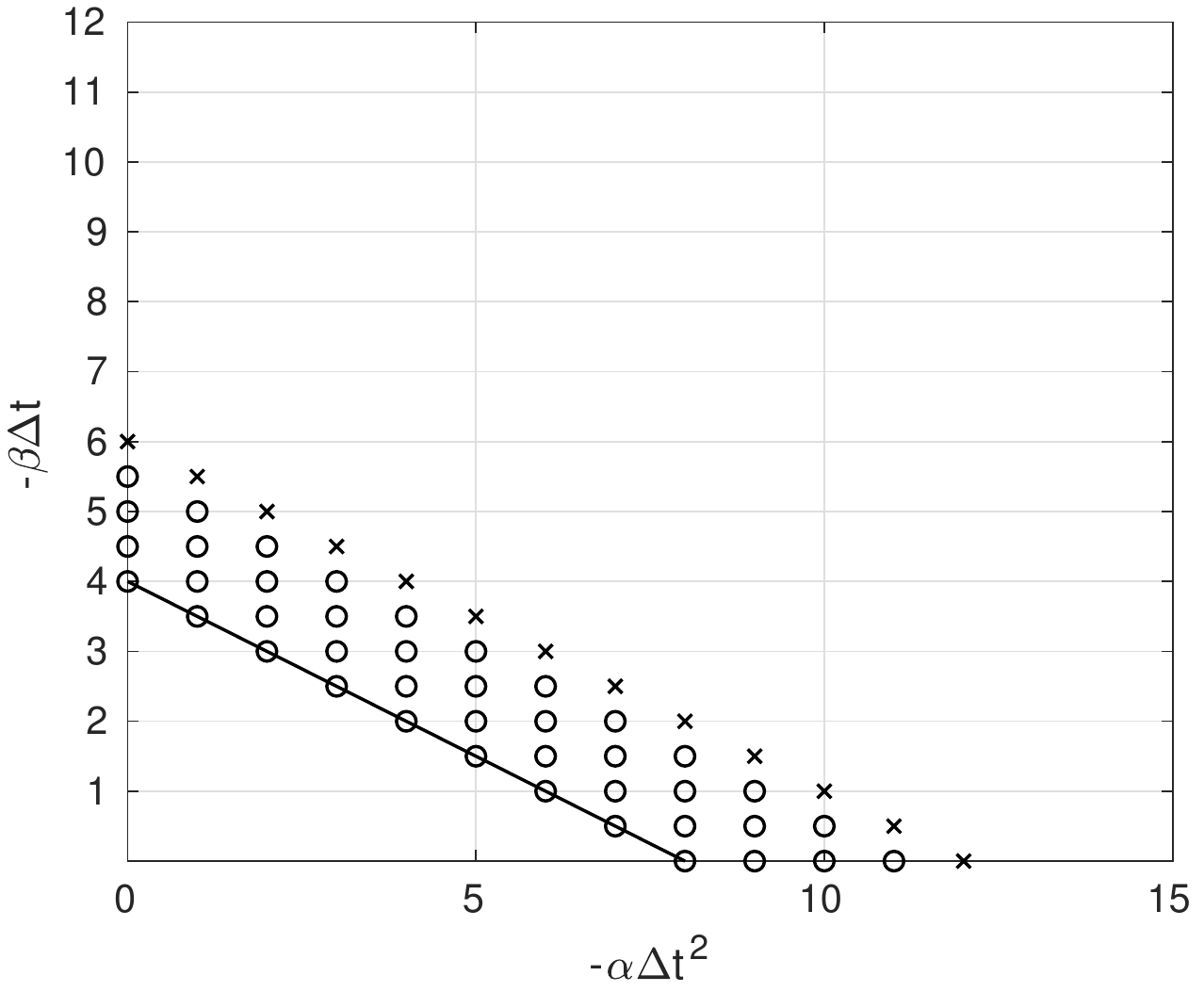}
       \put(50,84){\small{a)}}
      \end{overpic}
 \begin{overpic}[width=0.45\textwidth]{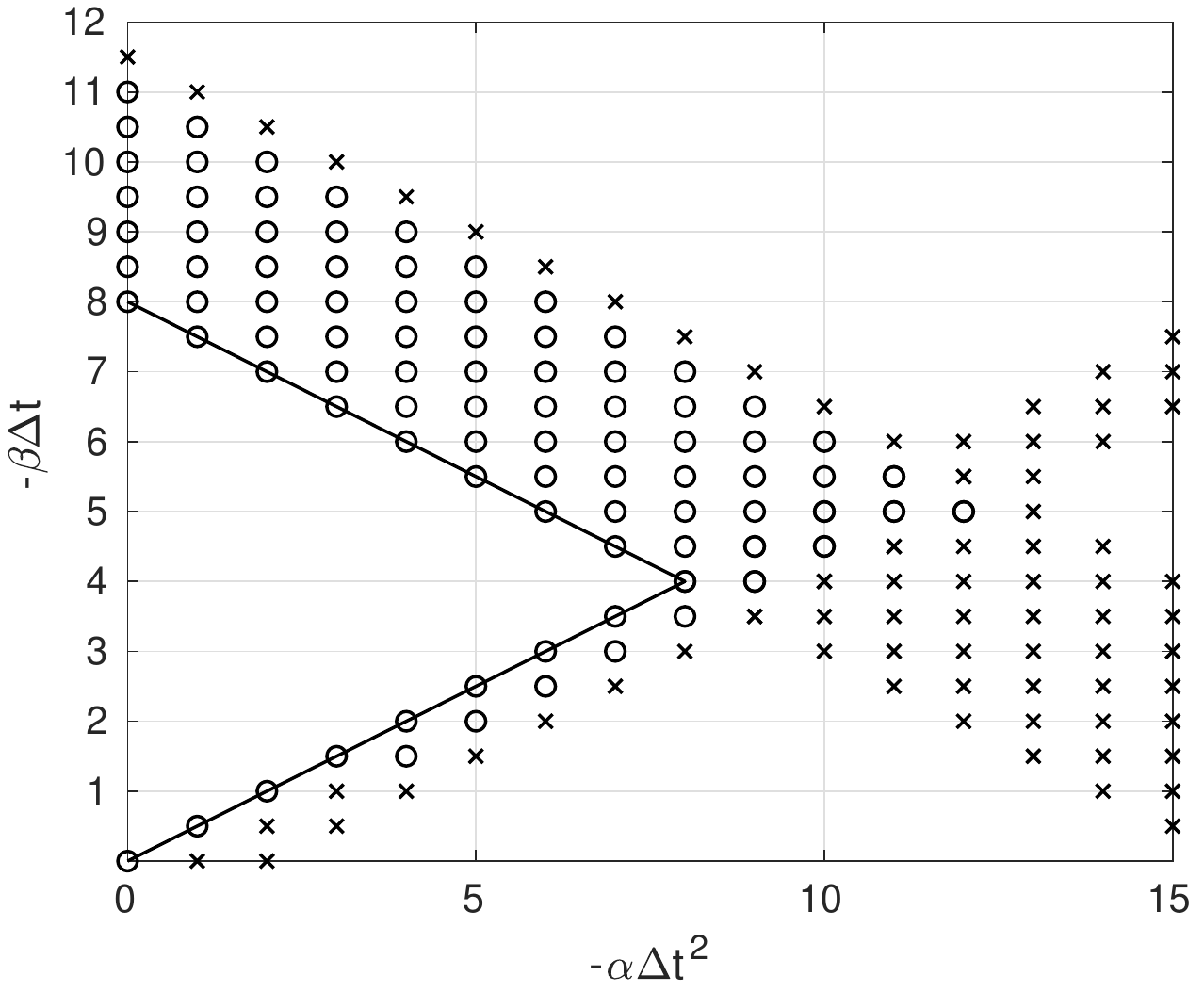}
       \put(50,84){\small{b)}}
   \end{overpic}
\caption{Numerical stability regions in 2D flow for $\gamma = 0$. Circles/crosses denote stable/unstable cases, respectively. Continuos lines denote the corresponding analytical curves: a) BDF1 time scheme, b) BDF2 time scheme.}
\label{fig:BDF_stability_2}
\end{figure}

\begin{figure}
\centering
 \begin{overpic}[width=0.45\textwidth]{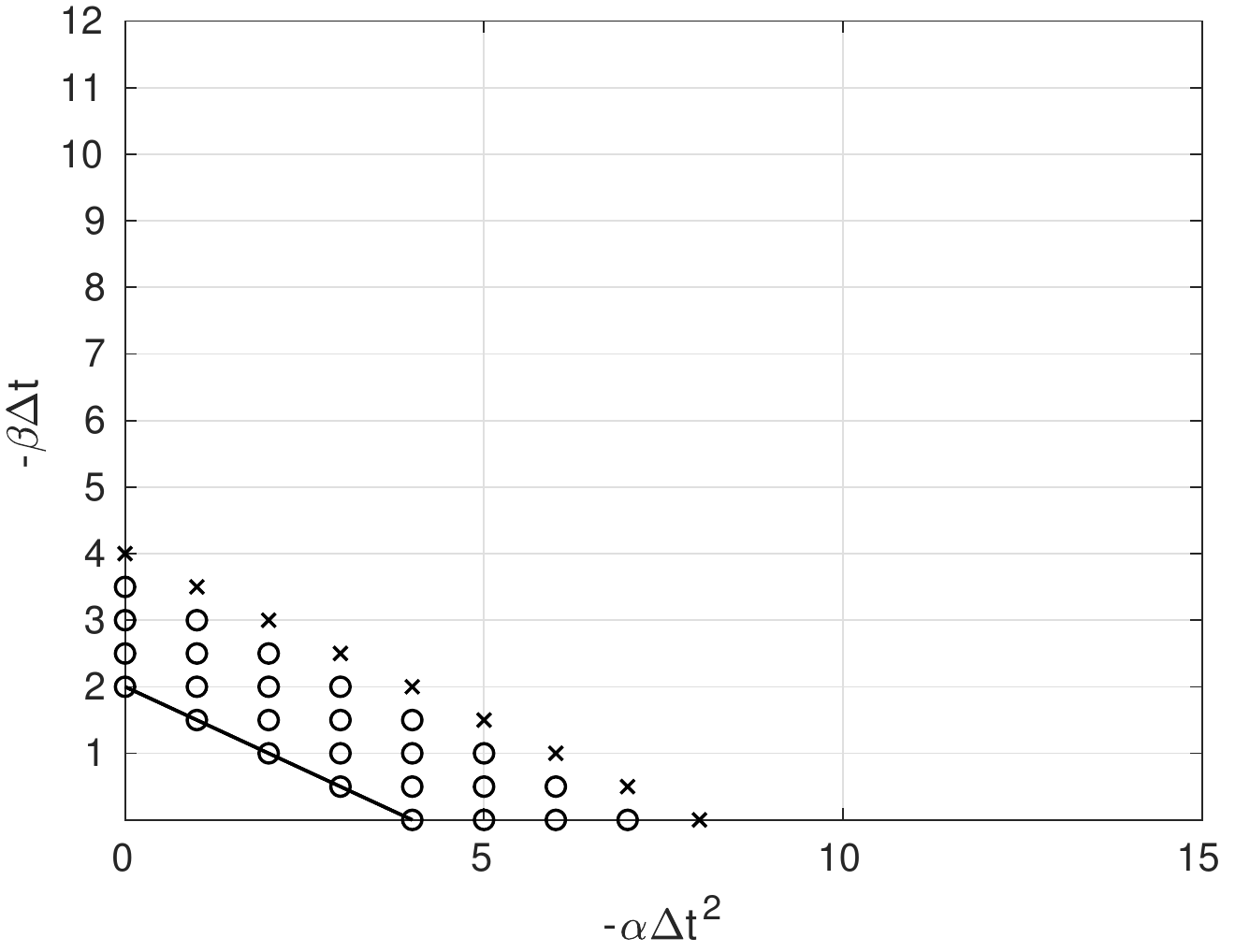}
        \put(50,82){\small{a)}}
      \end{overpic}
 \begin{overpic}[width=0.45\textwidth]{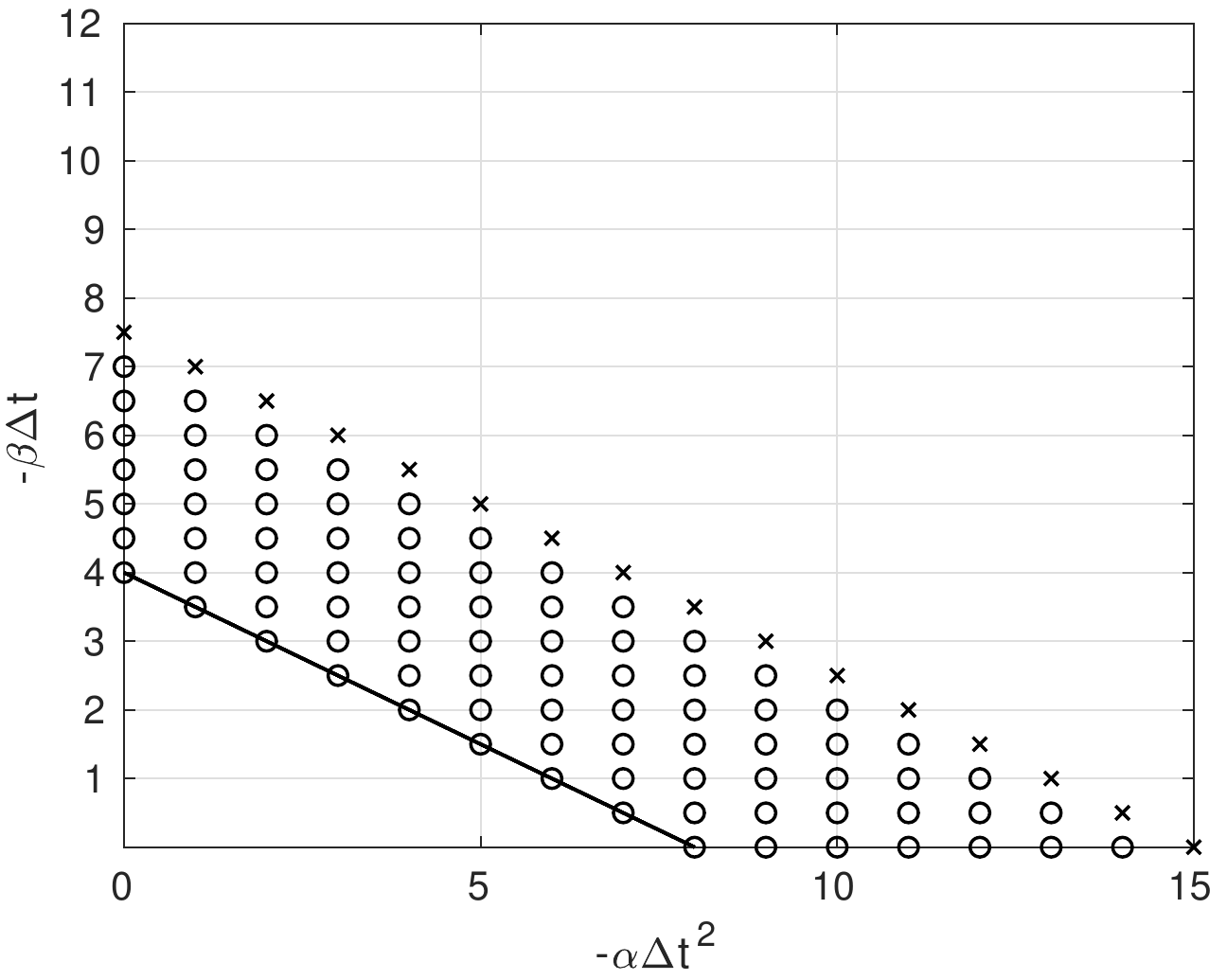}
        \put(50,82){\small{b)}}
   \end{overpic}
\caption{Numerical stability regions in 2D flow for $\gamma \neq 0$. Circles/crosses denote stable/unstable cases, respectively. Continuos lines denote the corresponding analytical curves: a) BDF1 time scheme and $-\gamma = 1$, b) BDF2 time scheme and $-\gamma = 2$.}
\label{fig:BDF_stability_3}
\end{figure}

Stability regimes \eqref{eq:stab1BDF2}-\eqref{eq:stab2BDF2} in 2D flow, obtained by set $C_{max} = \dfrac{1}{2}$ \cite{Shin2008},

\begin{equation}\label{eq:stab3BDF2}
\begin{cases}
- \alpha \Delta t^2 - 2 \beta \Delta t \leq 16 \text{    and    } \alpha \Delta t^2 - 2 \beta \Delta t \geq 0 & \text{for $-\gamma = 0$}, \\
- \alpha \Delta t^2 - 2 \beta \Delta t \leq 12.8 \text{    and    } \alpha \Delta t^2 - 14 \beta \Delta t \geq 0 & \text{for $-\gamma = 0.8$}, \\    
-\alpha \Delta t^2 - 2\beta \Delta t \leq 12 & \text{for $-\gamma = 1$}, \\
-\alpha \Delta t^2 - 2\beta \Delta t \leq 11.2 & \text{for $-\gamma = 1.2$}, \\
-\alpha \Delta t^2 - 2\beta \Delta t \leq 8 & \text{for $-\gamma = 2$}, \\
-\alpha \Delta t^2 - 2\beta \Delta t \leq 4 & \text{for $-\gamma = 3$}, \\
-\alpha \Delta t^2 = -\beta \Delta t = 0 & \text{for $-\gamma = 4$}, \\
\end{cases}
\end{equation}
are displayed in Figure \ref{fig:BDF_stability_1} b). The flow is stable in the region inside the line(s) and unstable in the region outside the line(s). Unlike the BDF1 time
scheme, for the BDF2 time scheme $\gamma$ affects not only the size of the stability region but also
its shape. We observe that for $\gamma = 0$ the stability region is limited by a superior line that has the same slope of that obtained for the BDF1 time scheme but greater intercept, and an inferior one that has equal and opposite slope with respect to the superior one and null intercept. When $-\gamma$ is approaching to $1$ from the left, we observe that the slope of the inferior line as well as the intercept of the superior one reduce. It is worth to point out that for $-\gamma < 1$ it is not possible to choose as forcing gains $\alpha \neq 0$ and $\beta = 0$ because no portion of the $-\alpha\Delta t^2$ axis belongs to the stability region. On the other hand, it is possible to choose values of $-\alpha \Delta t^2$ larger than ones related to the BDF1 time scheme by introducing low values of $-\beta\Delta t$ as well as very larger values of $-\beta\Delta t$ for medium-low values of $-\alpha\Delta t^2$. For $-\gamma = 1$, the inferior line concides with the portion of the $-\alpha\Delta t^2$ axis for which $-\alpha\Delta t^2 \leq 12$. Then, for $-\gamma > 1$, the shape of the stability region becomes the same of that one obtained for the BDF1 time scheme whilst the size is bigger than the biggest one related to the BDF1 time scheme that is obtained for $\gamma = 0$. The two regions coincide for $-\gamma = 2$. Finally, for $-\gamma > 2$, it continues to reduce until to degenerate into the axis origin for $-\gamma = 4$. Simulations of transverse oscillations of a circular cylinder are performed for $\gamma$ = 0 (Figure \ref{fig:BDF_stability_2} b)) and $-\gamma$ = 2 (Figure \ref{fig:BDF_stability_3} b)). Numerical stable and unstable cases are denoted by circles and crosses, respectively. Just like for the BDF1 time scheme, we obtain numerical stability regions wider than the analytical ones. By comparing Figure \ref{fig:BDF_stability_2} a) and Figure \ref{fig:BDF_stability_3} b), we observe that, thanks to the derivative action, the BDF2 time scheme allows to increase the maximum time step of about 13 \% with respect to the BDF1 time scheme for a given $\alpha$ when $\beta = 0$,
\begin{equation}\label{eq:BDF2vsBDF1}
\dfrac{\left(\Delta t\right)_{BDF2}}{\left(\Delta t\right)_{BDF1}} \approx \sqrt{\dfrac{14}{11}} \approx 1.13.
\end{equation}

\begin{figure}
\centering
\includegraphics[width=0.5\textwidth]{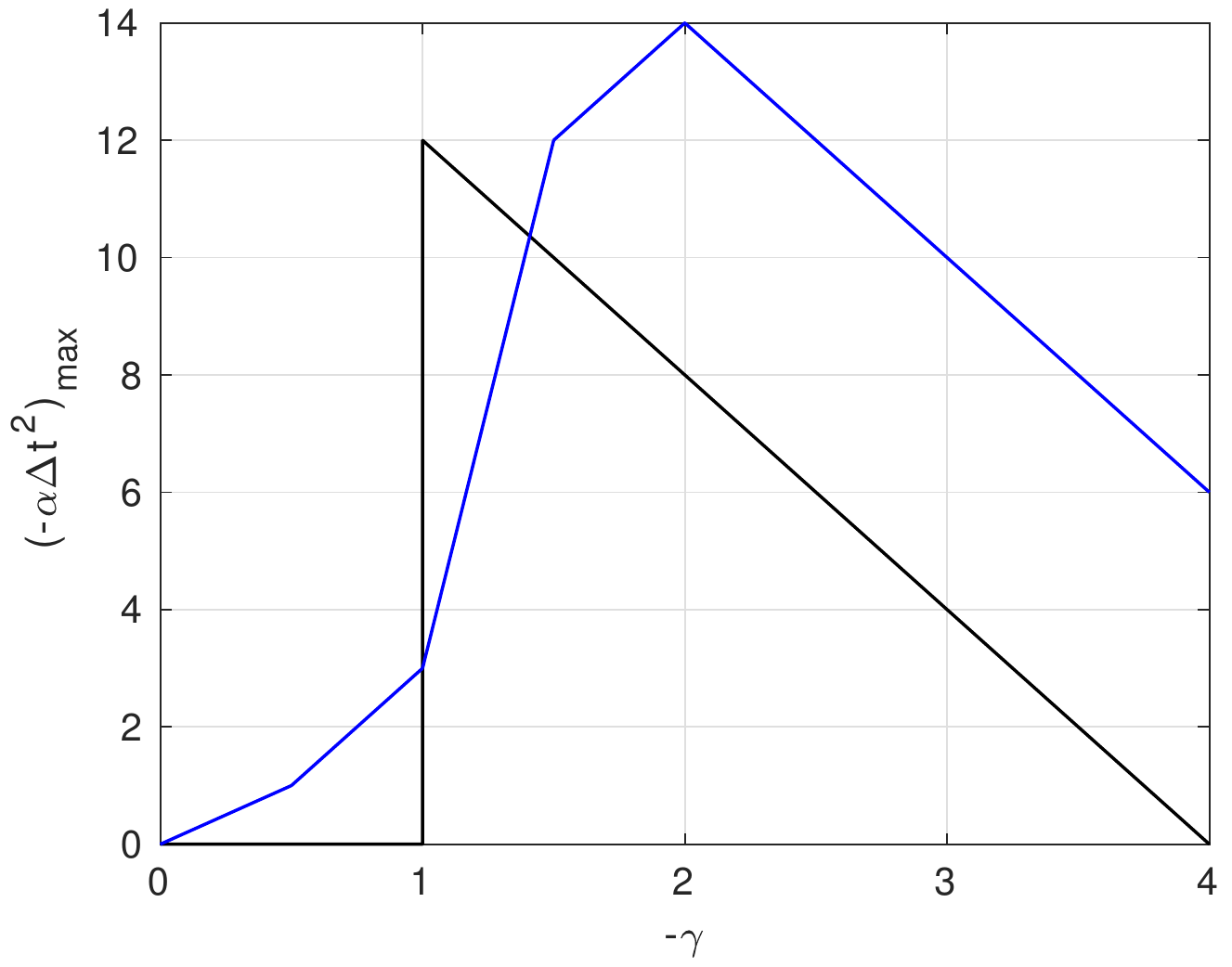}
\caption{Evolution of $\left(-\alpha \Delta t^2\right)_{max}$ over $-\gamma$ for 2D flow related to the BDF2 time scheme: comparison between analytical (black line) and numerical (blue line) curves.}
\label{fig:alpha_max}
\end{figure}

Of course, just like for the BDF1 scheme, the marginal stability line is expected to change depending on the kind of flow. Now, we investigate in detail how the numerical system behaves around $-\gamma = 1$ because, as previously noted, for this value the BDF2 time scheme degenerates into the BDF1 time scheme unless a scale factor, and the shape of the stability region changes. In Figure \ref{fig:alpha_max} numerical and analytical evolutions of the maximum value of $-\alpha \Delta t^2$ over $-\gamma$ for which the stability occurs are displayed. We can observe the discontinuity exhibited by the analytical curve for $-\gamma = 1$, associated to the changing of shape of the stability region. Numerically, the behaviour is smoother and we observe an early onset of instability at values of $-\alpha \Delta t^2$ lower than analytical predictions for $-\gamma$ next to 1 from right. 

In conclusion, we learned that the when the BDF2 time scheme is used the PID controller with an \emph{ad hoc} choice of $\gamma$ allows to obtain larger stability regions with respect to the combination of the standard PI controller with the BDF1 scheme. Therefore, the derivative action provides a general improvement of the stability features of the system when higher order time advancing schemes are considered.

\section{Numerical results}\label{sec:numerical_results}

In this section, we present several numerical results for the VBM. Four 2D benchmarks, involving stationary or ridigly moving bodies, are investigated: flow past a stationary cylinder, oscillatory flow past a stationary cylinder, transverse oscillation of a cylinder in a free-stream and inline oscillation of a cylinder in a fluid at rest. Finally, we will present a FSI benchmark related to the natural oscillation of a submerged cantilever rectangular beam in a fluid at rest. The number of PISO \cite{PISO} loops has been fixed to 2 for all the simulations. The linear algebraic system associated with the momentum equation \eqref{eq:evolveFV-1.1} is solved using an iterative solver with symmetric Gauss-Seidel smoother. For Poisson problem \eqref{eq:Poisson} we use Geometric Agglomerated Algebraic Multigrid Solver GAMG with the Gauss-Seidel smoother. The required accuracy is 1e-7 at each time step.

\begin{table}[h]
\centering
\begin{tabular}{ccccc}
\multicolumn{2}{c}{} \\
\cline{1-5}
mesh name &  Domain
& $h_{min}$ & 
$h_{max}$ & No. of cells  \\
\hline
$117\mbox{k}$ & $[-50 \quad 50]^2$ & 1.63e-2 & 9.76 & 117,312   \\
$260\mbox{k}$ & $[0 \quad 8]^2$ & 1.56e-2 & 1.56e-2 & 262,144  \\
$1000\mbox{k}$ & $[0 \quad 16]^2$ & 1.56e-2 & 1.56e-2 & 1,048,576 \\
\hline
\end{tabular}
\caption{Name, size of the computational domain, minimum diameter $h_{min}$, maximum diameter $h_{max}$, and number of cells for all the meshes used for the tests presented in Secs. \ref{sec:staz_cyl}, \ref{sec:osc_y_cyl} and \ref{sec:osc_x_cyl}.}
\label{tab:mesh}
\end{table}



\begin{table}[h]
\centering
\tiny
\begin{tabular}{l
ccc
ccccc
cccccc}
\multicolumn{2}{c}{} \\
\cline{1-12}

& & Mesh & $-\alpha$ & $-\beta$&  $-\gamma$ & $\bar{C}_d$ & $C'_l$ & $\mbox{St}$  & $\Delta t$ & $CFL$ & Time scheme  \\
\hline
Case 1 & present & $260\mbox{k}$ & 4.8e4 & 0 & 0  & 1.57   & 0.44     & 0.159  & 1.2e-2  & 1.35 & BDF1      \\
 &  present & $260\mbox{k}$ & 4.8e4 & 0 & 2  & 1.6   & 0.52     & 0.163  & 1.2e-2  & 1.35 & BDF2      \\
 & \cite{Shin2008} & $260\mbox{k}$ & 4.8e4 & 0 & 0 & 1.44  & 0.35  & 0.168 & 1.2e-2  & 1.35 & BDF1      \\ \hline
Case 2  & present & $260\mbox{k}$ & 4.8e4 & 0 & 0   &  1.56  &  0.45    & 0.162 & 6.0e-3    & 0.7 & BDF1      \\
  &  present &$260\mbox{k}$ & 4.8e4 & 0 & 2   & 1.58  &  0.49    & 0.164  & 6.0e-3    & 0.7 & BDF2      \\
 & \cite{Shin2008}  &  $260\mbox{k}$ & 4.8e4 & 0  & 0 & 1.44  & 0.35  & 0.168 & 6.0e-3 & 0.7 & BDF1      \\\hline
Case 3  & present & $1000\mbox{k}$ &  4.8e4 & 0 & 0   &  1.56  &  0.45    & 0.162 & 6.0e-3    & 0.7 & BDF1      \\
 &  present & $1000\mbox{k}$ & 4.8e4 & 0 & 2   & 1.57  &  0.48    & 0.165 & 6.0e-3    & 0.7 & BDF2      \\
 & \cite{Shin2008}  & $1000\mbox{k}$ & 4.8e4 & 0 & 0 & 1.37  & 0.34  & 0.163  & 6.0e-3 & 0.7 & BDF1      \\\hline
Case 4 & present & $117\mbox{k}$ & 4.8e4 & 0 &  0  &  1.34  & 0.33   & 0.16 &  1.2e-2  & 1.15   & BDF1      \\
   & present & $117\mbox{k}$ & 4.8e4 & 0 &  2  &  1.35  & 0.35   & 0.161  & 1.2e-2  & 1.15    & BDF2      \\
   & present & $117\mbox{k}$ & 4.8e4 & 0 &  0  &  1.34  & 0.33   & 0.16  & 1.6e-2  & 1.55    & BDF1      \\
   & present & $117\mbox{k}$ & 4.8e4 & 0 &  2  &  1.36  & 0.36  &  0.161 & 1.8e-2  & 1.75    & BDF2      \\\hline
   & \cite{Lai2000} & $260\mbox{k}$  & 4.8e4 & 0 & 0 & 1.52 & 0.29 & 0.155  & 1.8e-3  &   &  Crank–Nicolson     \\
   & \cite{Uhlmann2005} &   &    &    &     &   1.50   & 0.35     &  0.172 & 3e-3    \\
   & \cite{Kim2001} &   &    &    &     &  1.33    & 0.32     & 0.165      \\ 
   & \cite{Constant2017} &   &    &    &     &  1.38    &      & 0.165 & & 0.5 & BDF1     \\
   & \cite{Lee2006}      &   &    &    &     &  1.33    & 0.28  & 0.166      \\ 
   & \cite{Linnick2005} &   &    &    &     &   1.34   & 0.33     &  0.166    \\
   & \cite{Huang2007_2} &   &    &    &     &   1.36   & 0.33     &  0.167    \\
\hline
\end{tabular}
\caption{Flow past a stationary cylinder at $\mbox{Re} = 100$: values of aerodynamic coefficients, $C'_l$ and $\bar{C}_d$, and Strouhal number, $\mbox{St}$, compared against the results reported in other studies.}
\label{tab:cyl_staz_1}
\end{table}

\begin{table}[h]
\centering
\begin{tabular}{l
ccc
ccccc
c}
\multicolumn{2}{c}{} \\
\cline{1-4}

$-\alpha\Delta t^2$  & $\bar{C}_d$ & $C'_l$ & $\mbox{St}$ \\
\hline
0     &  1.4151   &   0.0532   & 0.1493       \\
0.01  &  1.3376   &   0.2715   & 0.1672       \\
0.1   &  1.3147   &   0.2939   & 0.1599      \\
0.5   &  1.3348   &   0.3274   & 0.1601      \\
1     &  1.3375   &   0.3311   & 0.16      \\
3.5   &  1.3389   &   0.3333   & 0.16     \\
\hline
\end{tabular}
\caption{Flow past a stationary cylinder at $\mbox{Re} = 100$: influence of $-\alpha \Delta t^2$ on values of aerodynamic coefficients, $C'_l$ and $\bar{C}_d$, and Strouhal number, for $-\beta \Delta t= 1.5$ and $-\gamma = 0$, and the BDF1 time scheme.
}
\label{tab:cyl_staz_alpha}
\end{table}

\begin{table}[h]
\centering
\begin{tabular}{l
ccc
ccccc
c}
\multicolumn{2}{c}{} \\
\cline{1-4}

 $-\beta\Delta t$  & $\bar{C}_d$ & $C'_l$ & $\mbox{St}$   \\
\hline
 0   &  1.3215   &   0.3088   & 0.1602      \\
 0.5    &  1.3219   &   0.3087   & 0.1603       \\
 1  &    1.3221   &   0.3085   & 0.1604     \\
 1.5   &  1.3147   &   0.2939   & 0.1599      \\
 2    &  1.3231   &   0.3092   & 0.1606     \\
 2.5   &  1.3239   &   0.3099   & 0.1607       \\
 3.5    &  1.3253   &   0.3110   & 0.1608       \\
\hline
\end{tabular}
\caption{Flow past a stationary cylinder at $\mbox{Re} = 100$: influence of $-\beta \Delta t$ on values of aerodynamic coefficients, $C'_l$ and $\bar{C}_d$, and Strouhal number, for $-\alpha \Delta t^2= 0.1$ and $-\gamma = 0$, and the BDF1 time scheme.
}
\label{tab:cyl_staz_beta}
\end{table}

\begin{table}[h]
\centering
\begin{tabular}{l
ccc
ccccc
c}
\multicolumn{2}{c}{} \\
\cline{1-4}
$-\gamma$ & $\bar{C}_d$ & $C'_l$ & $\mbox{St}$ \\
\hline
 0  &  1.3147   &   0.2939   & 0.1599     \\
 0.5  &  1.3156   &   0.2945   & 0.16     \\
 1  &  1.3165   &   0.301   & 0.1603       \\
\hline
\end{tabular}
\caption{Flow past a stationary cylinder at $\mbox{Re} = 100$: influence of $-\gamma$ on values of aerodynamic coefficients, $C'_l$ and $\bar{C}_d$, and Strouhal number, for $-\alpha \Delta t^2= 0.1$ and $-\beta \Delta t = 1.5$, and the BDF1 time scheme.}
\label{tab:cyl_staz_gamma}
\end{table}

\begin{table}[h]
\centering
\begin{tabular}{l
ccc
ccccc
c}
\multicolumn{2}{c}{} \\
\cline{1-5}

$\Delta s/h$ & $N_L$ & $\bar{C}_d$ & $C'_l$ & $\mbox{St}$  \\
\hline
0.1  &  1880 &  1.321   &  0.3271   &  0.1587    \\
0.25 & 770  &  1.324  &  0.328   &  0.1587      \\
0.5 &  385  &  1.3268  &  0.3304   &  0.1587    \\
1   &  194 &   1.3267   &  0.3351  &  0.1613     \\
2   &  97 &   1.3031   &  0.3156  &  0.1613      \\
\hline
\end{tabular}
\caption{Flow past a stationary cylinder at Re = 100: influence of the ratio of Lagrangian point distance to Eulerian grid width, $\Delta s/h$, on values of aerodynamic coefficients, $C'_l$ and $\bar{C}_d$, and Strouhal number, $\mbox{St}$, for $-\alpha \Delta t^2= 3.9$ and $-\beta \Delta t= 1.9$, and the BDF1 time scheme.}
\label{tab:cyl_staz_ds_h}
\end{table}

\subsection{Flow past a stationary cylinder}\label{sec:staz_cyl}
\begin{figure}
\centering
 \begin{overpic}[width=0.45\textwidth]{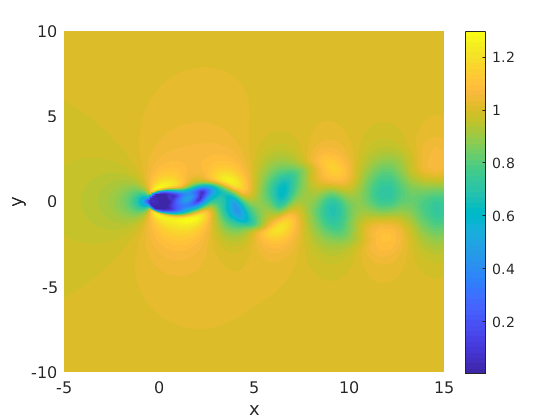}
       \put(50,78){\small{a)}}
      \end{overpic}
 \begin{overpic}[width=0.45\textwidth]{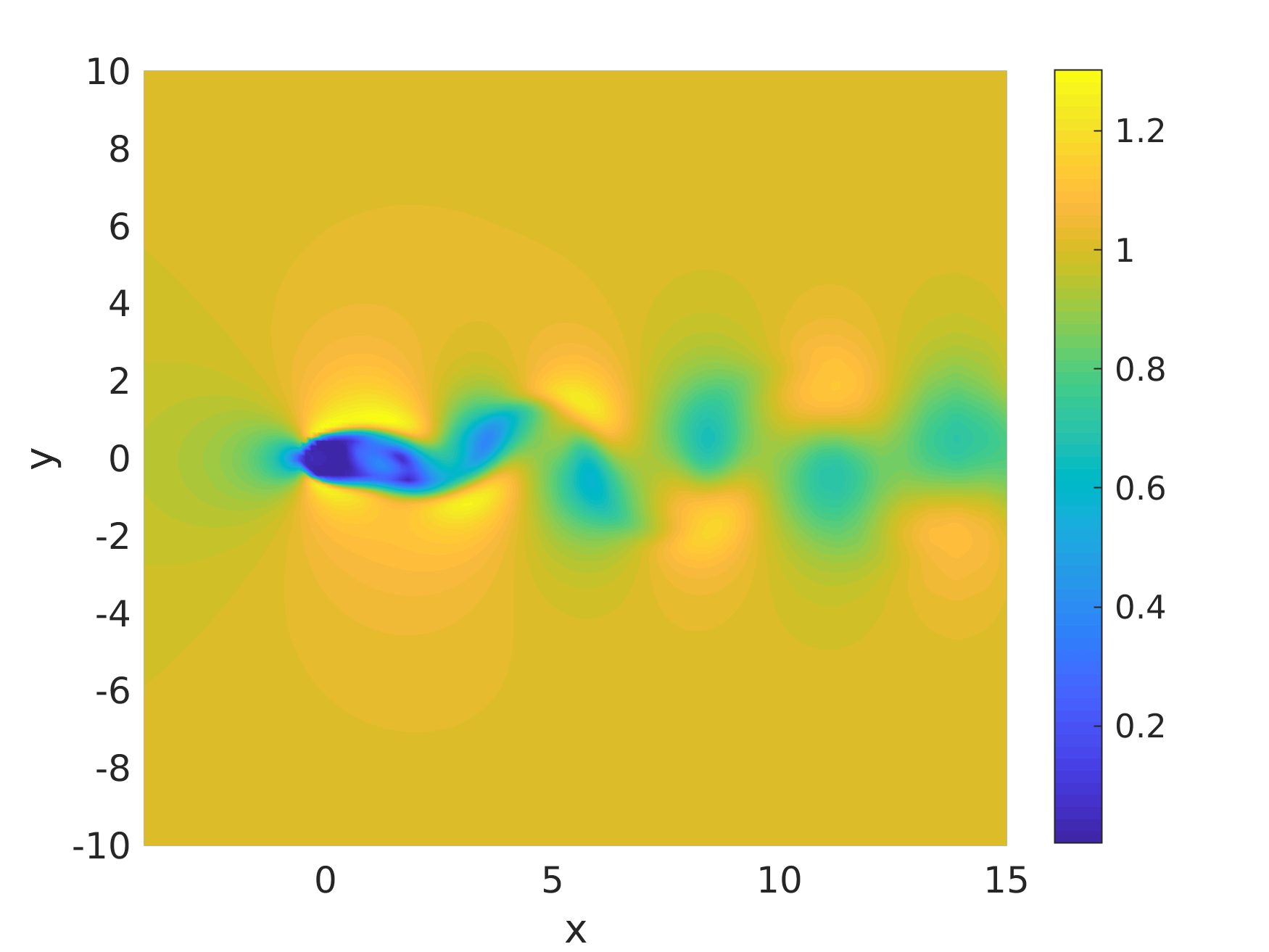}
       \put(50,78){\small{b)}}
   \end{overpic}
\caption{Flow past a stationary cylinder at Re = 100. Velocity magnitude for \emph{Case 4}: a) BDF1 time scheme, b) BDF2 time scheme.}
\label{fig:vel_mag}
\end{figure}

\begin{figure}
\centering
 \begin{overpic}[width=0.46\textwidth]{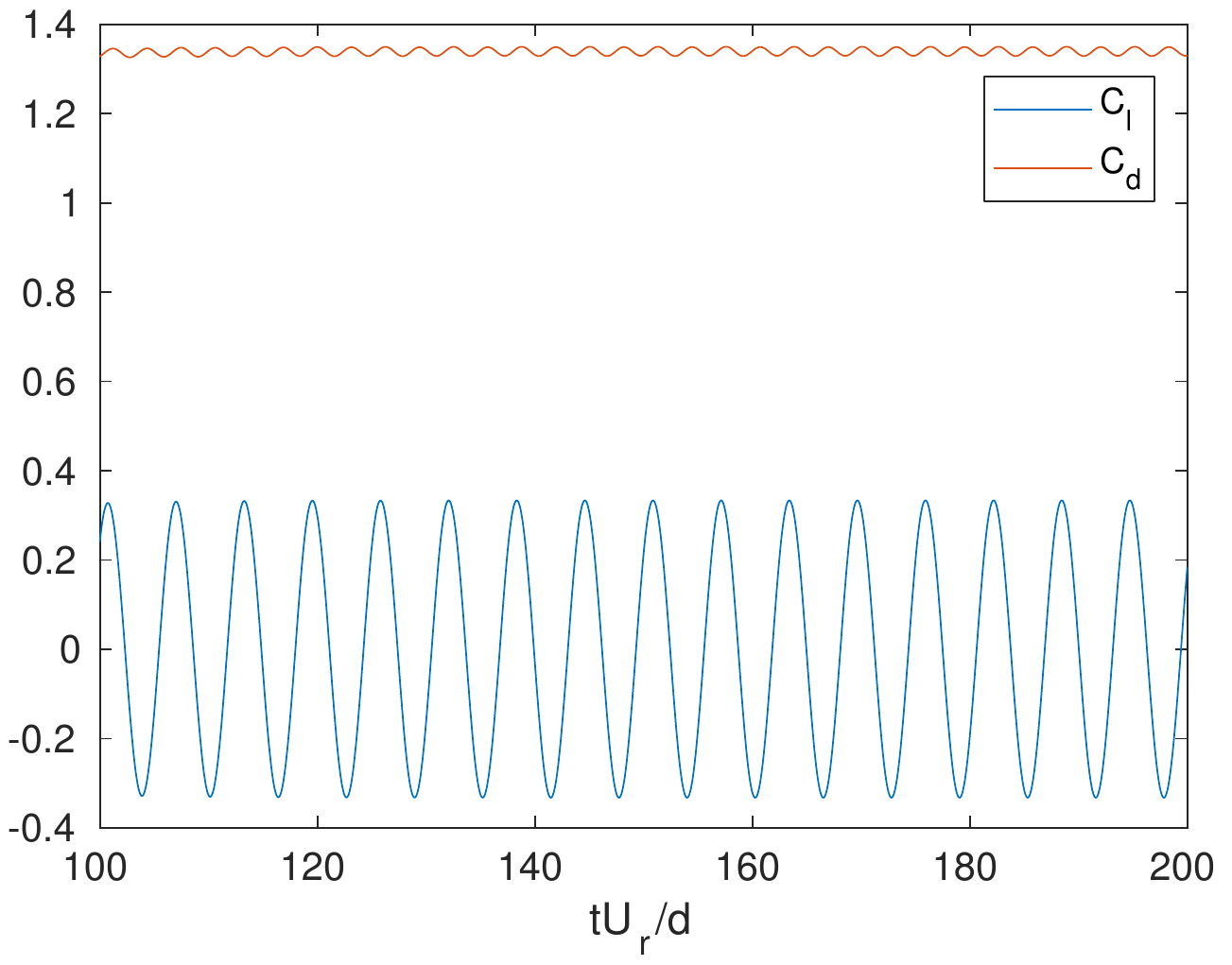}
       \put(50,84){\small{a)}}
      \end{overpic}
 \begin{overpic}[width=0.45\textwidth]{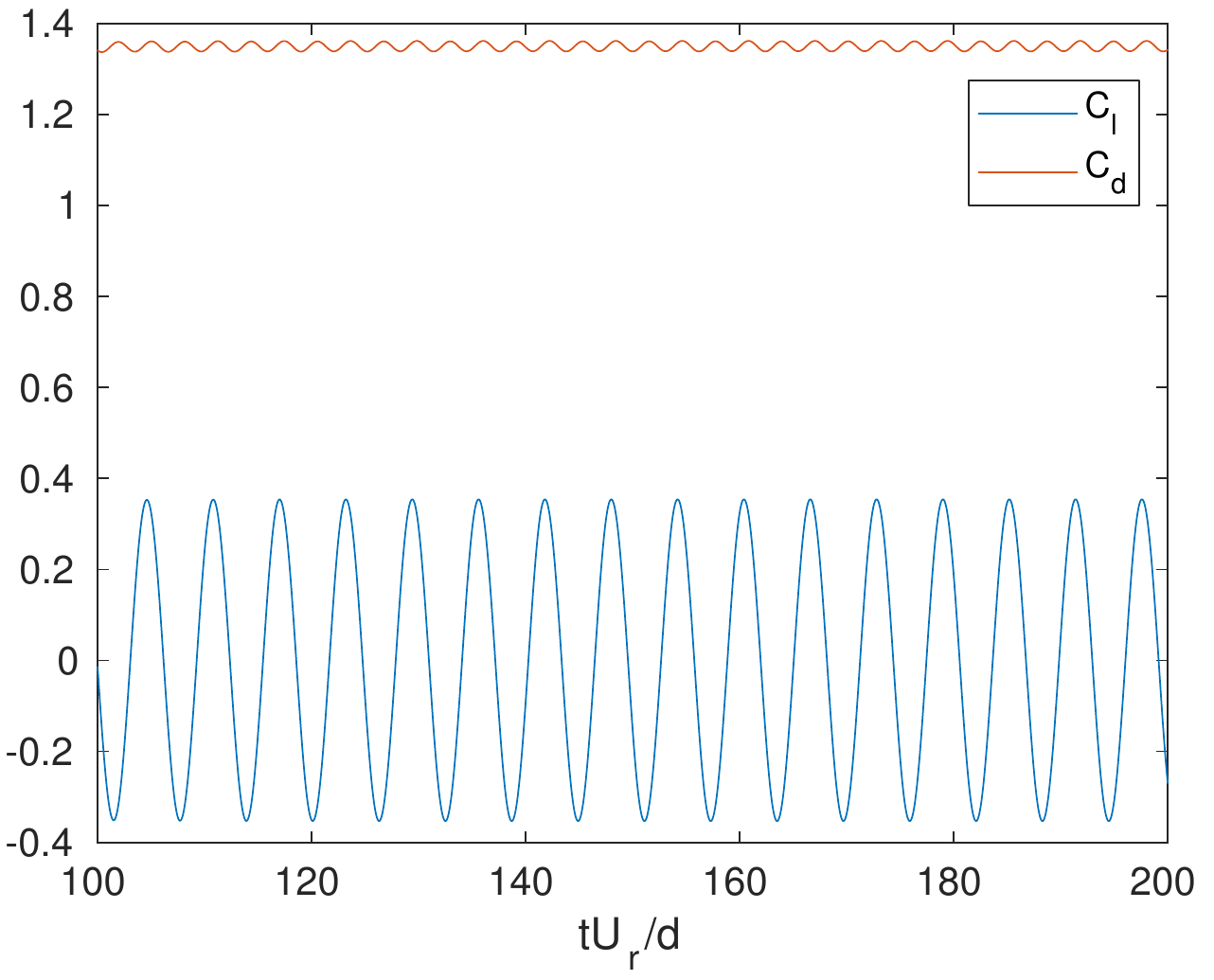}
       \put(50,85){\small{b)}}
   \end{overpic}
\caption{Flow past a stationary cylinder at Re = 100. Time history of drag and lift coefficients for \emph{Case 4}: a) BDF1 time scheme, b) BDF2 time scheme.}
\label{fig:coeffs}
\end{figure}

\begin{figure}
\centering
 \begin{overpic}[width=0.45\textwidth]{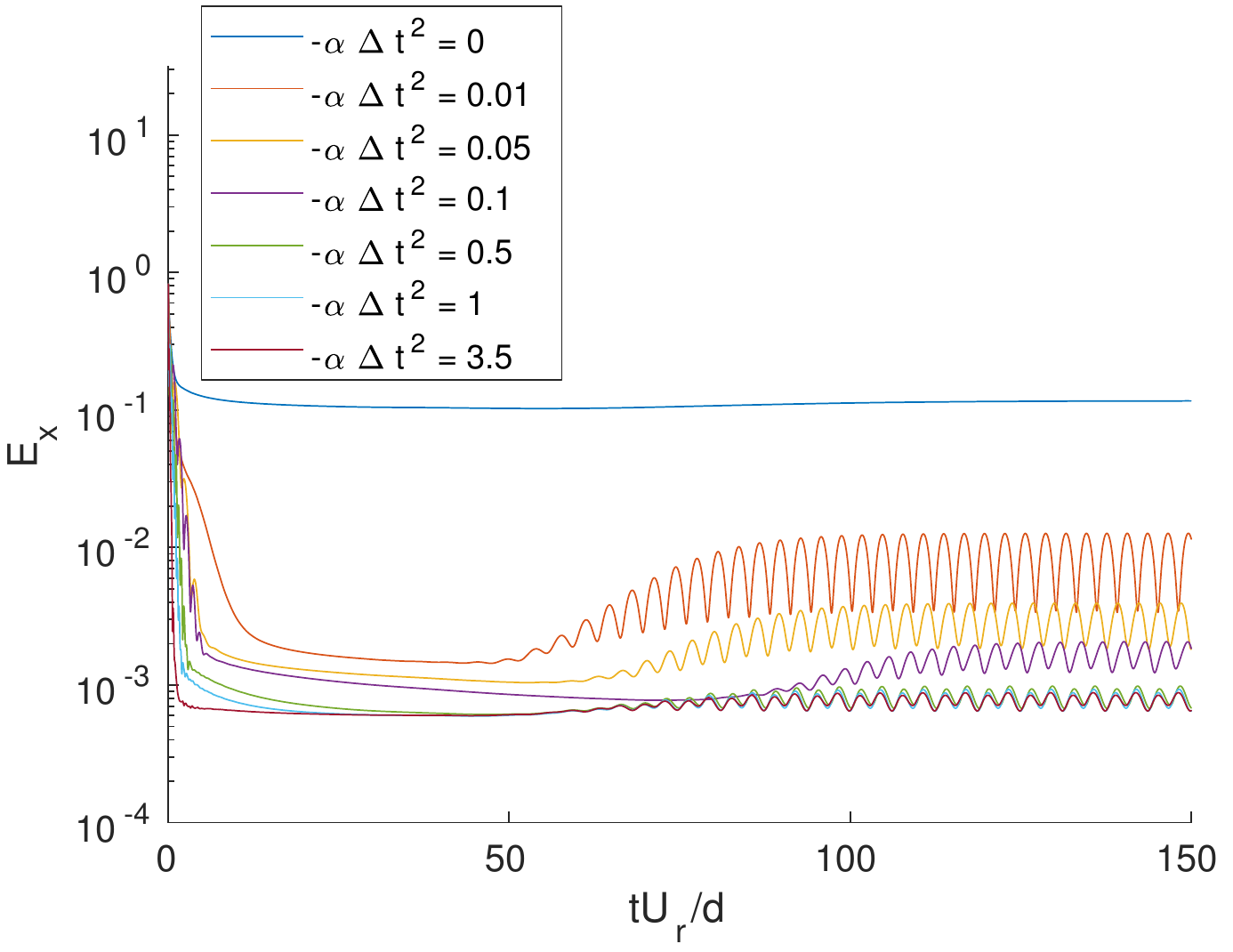}
       \put(32,78.5){\small{BDF1, -$\beta\Delta t$ = 1.5, $\gamma$ = 0}}
      \end{overpic}
 \begin{overpic}[width=0.45\textwidth]{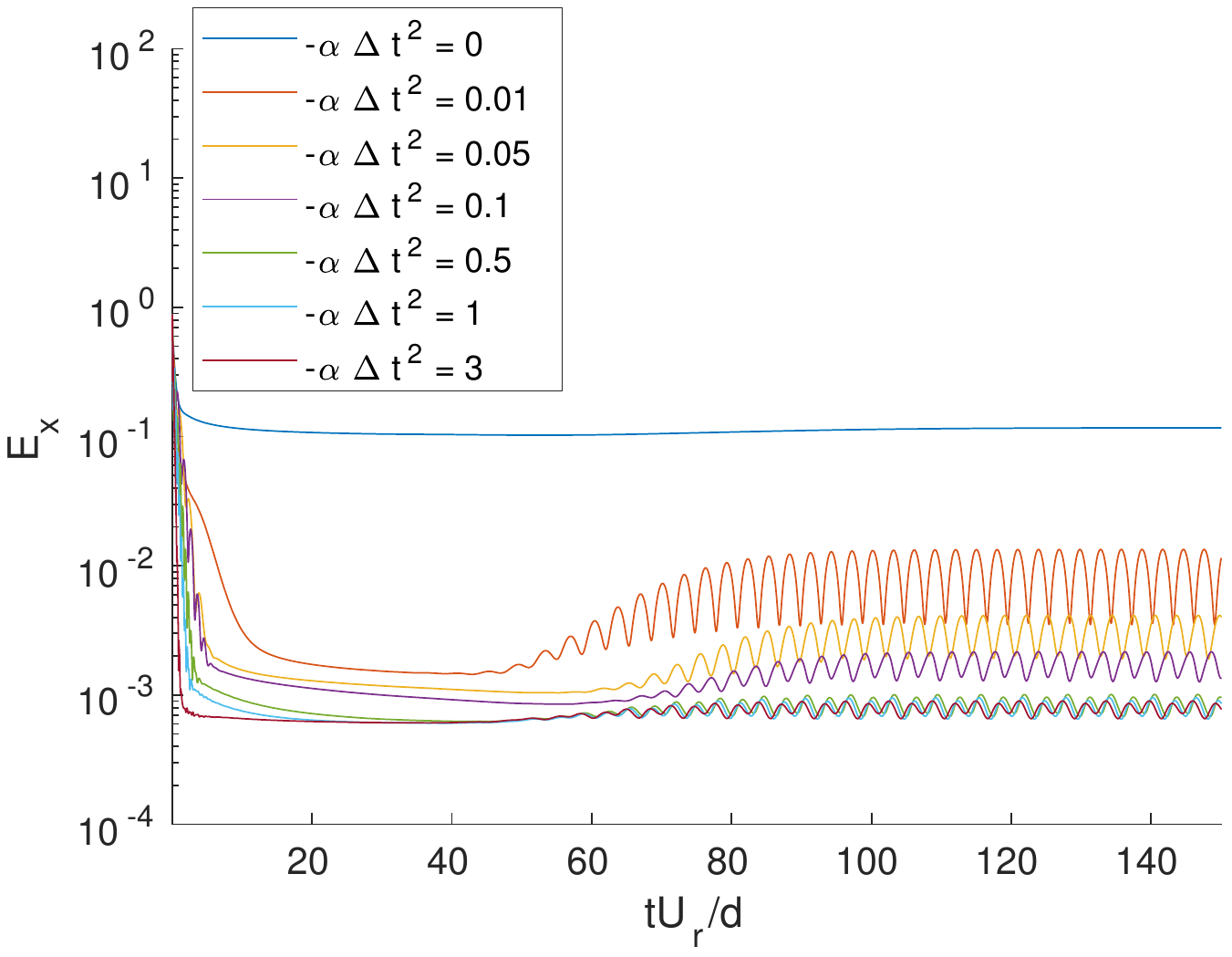}
       \put(32,78.5){\small{BDF2, -$\beta\Delta t$ = 1.5, $\gamma$ = 0}}
   \end{overpic}\\
\caption{Flow past a stationary cylinder at Re = 100: influence of $-\alpha \Delta t^2$ on the time evolution of $E_x$.}
\label{fig:Err_x_alpha}
\end{figure}

\begin{figure}
\centering
 \begin{overpic}[width=0.45\textwidth]{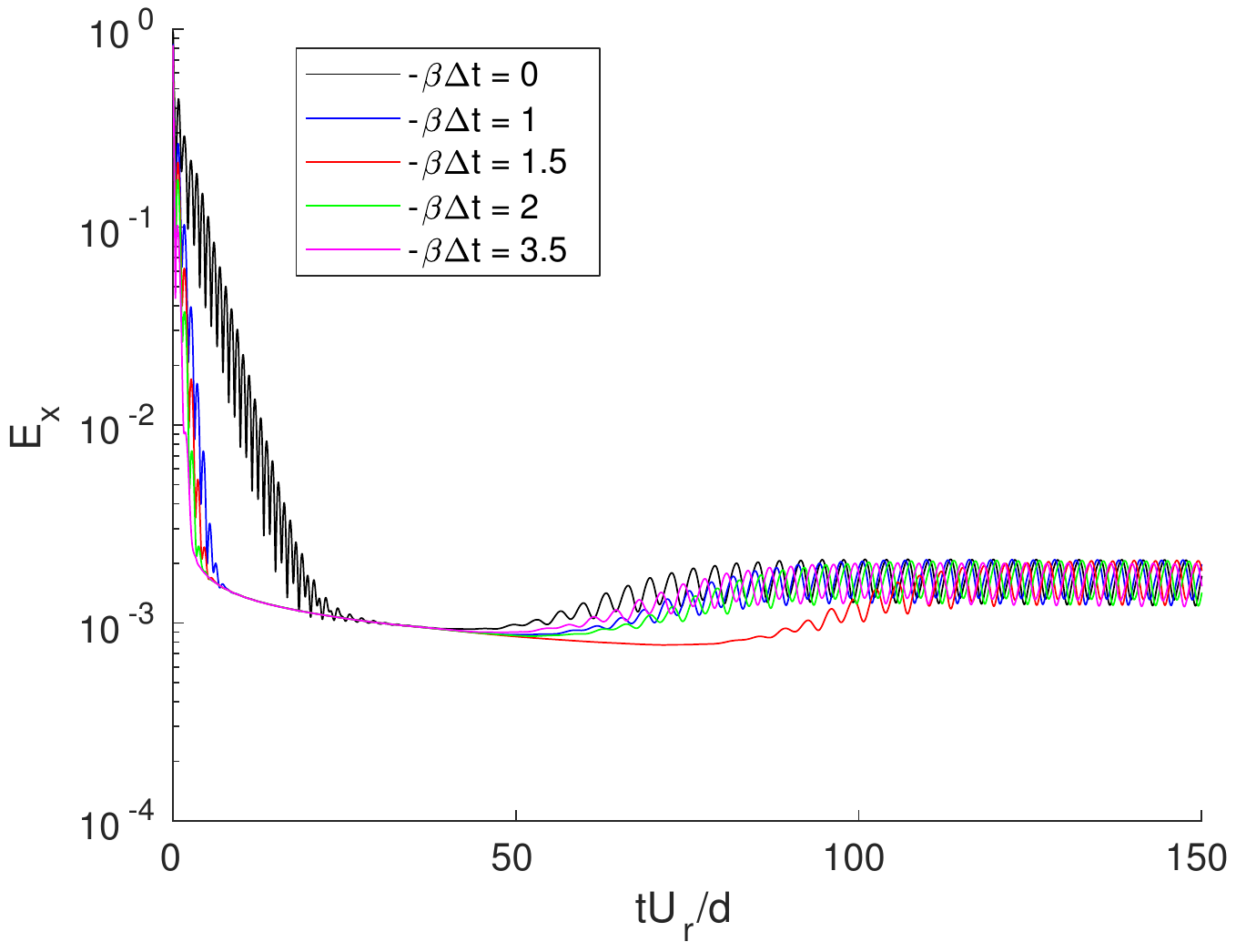}
       \put(32,75){\small{BDF1, -$\alpha\Delta t^2$ = 0.1, $\gamma$ = 0}}
      \end{overpic}
 \begin{overpic}[width=0.45\textwidth]{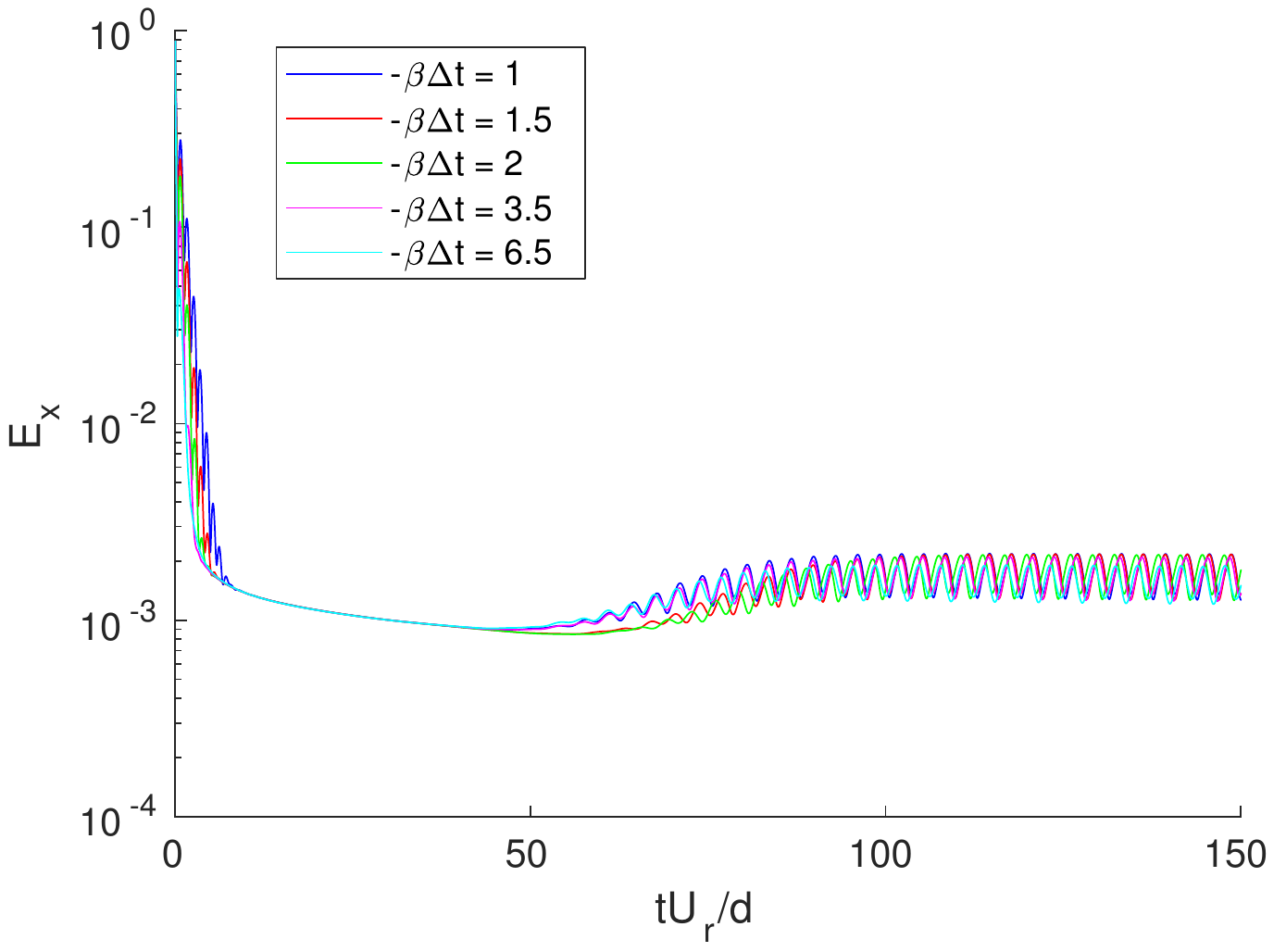}
       \put(32,75){\small{BDF2, -$\alpha\Delta t^2$ = 0.1, $\gamma$ = 0}}
   \end{overpic}\\
\caption{Flow past a stationary cylinder at Re = 100: influence of $-\beta \Delta t$ on the time evolution of $E_x$.}
\label{fig:Err_x_beta}
\end{figure}

\begin{figure}
\centering
 \begin{overpic}[width=0.45\textwidth]{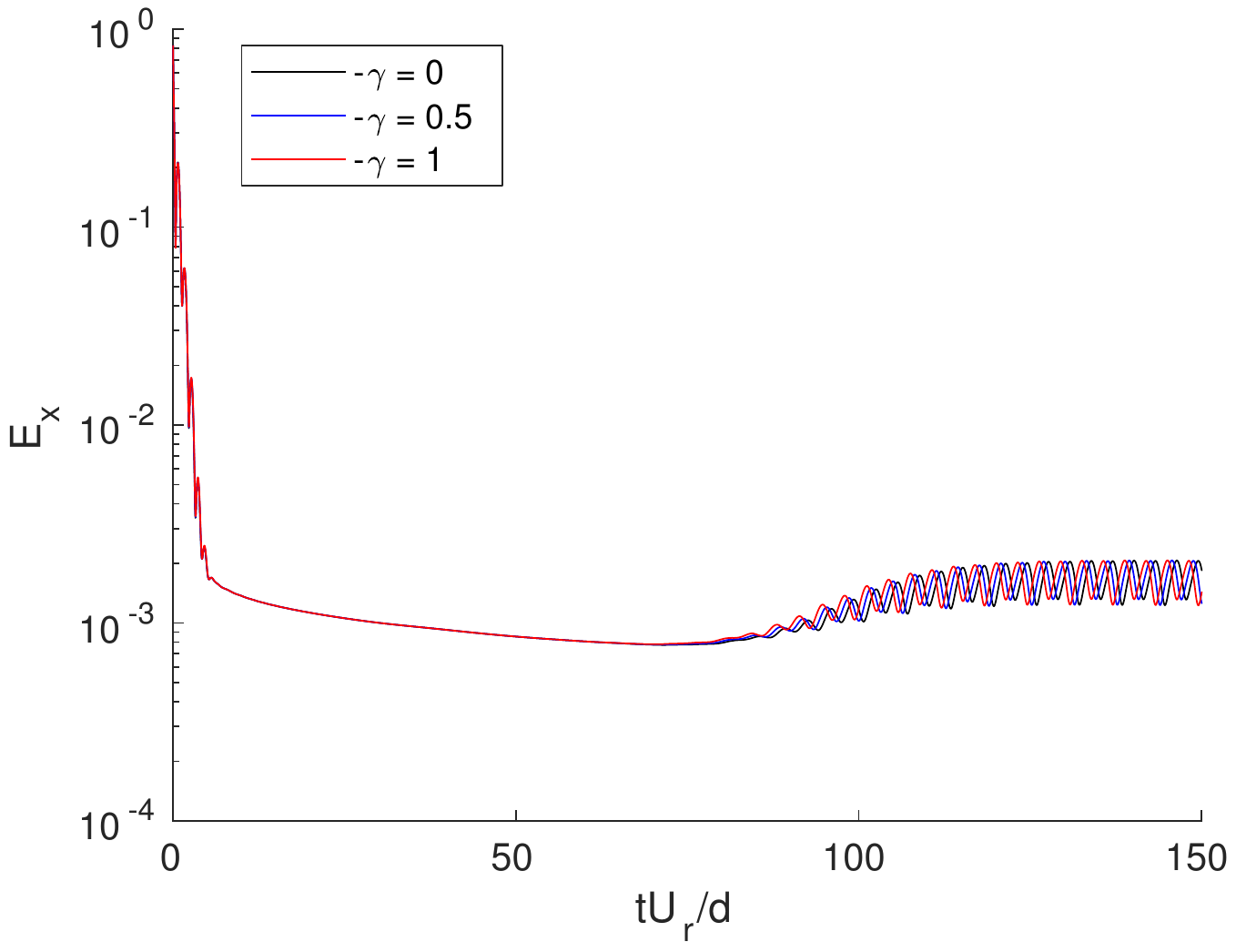}
       \put(27,75){\small{BDF1, -$\alpha\Delta t^2$ = 0.1, -$\beta\Delta t$ = 1.5}}
      \end{overpic}
 \begin{overpic}[width=0.45\textwidth]{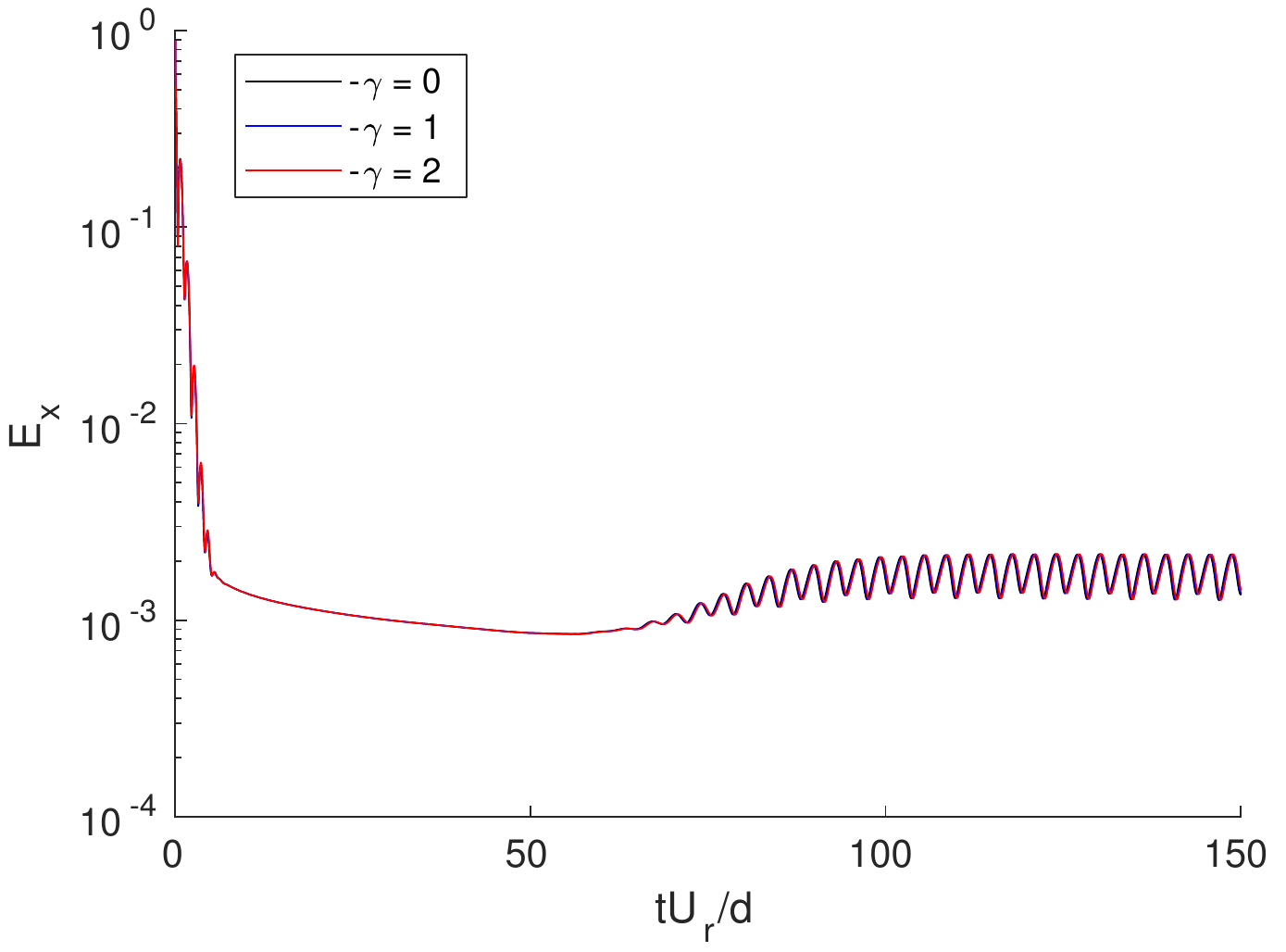}
       \put(27,75){\small{BDF2, -$\alpha\Delta t^2$ = 0.1, -$\beta\Delta t$ = 1.5}}
   \end{overpic}\\
\caption{Flow past a stationary cylinder at Re = 100: influence of $-\gamma$ on the time evolution of $E_x$.}
\label{fig:Err_x_gamma}
\end{figure}

\begin{figure}
\centering
 \begin{overpic}[width=0.45\textwidth]{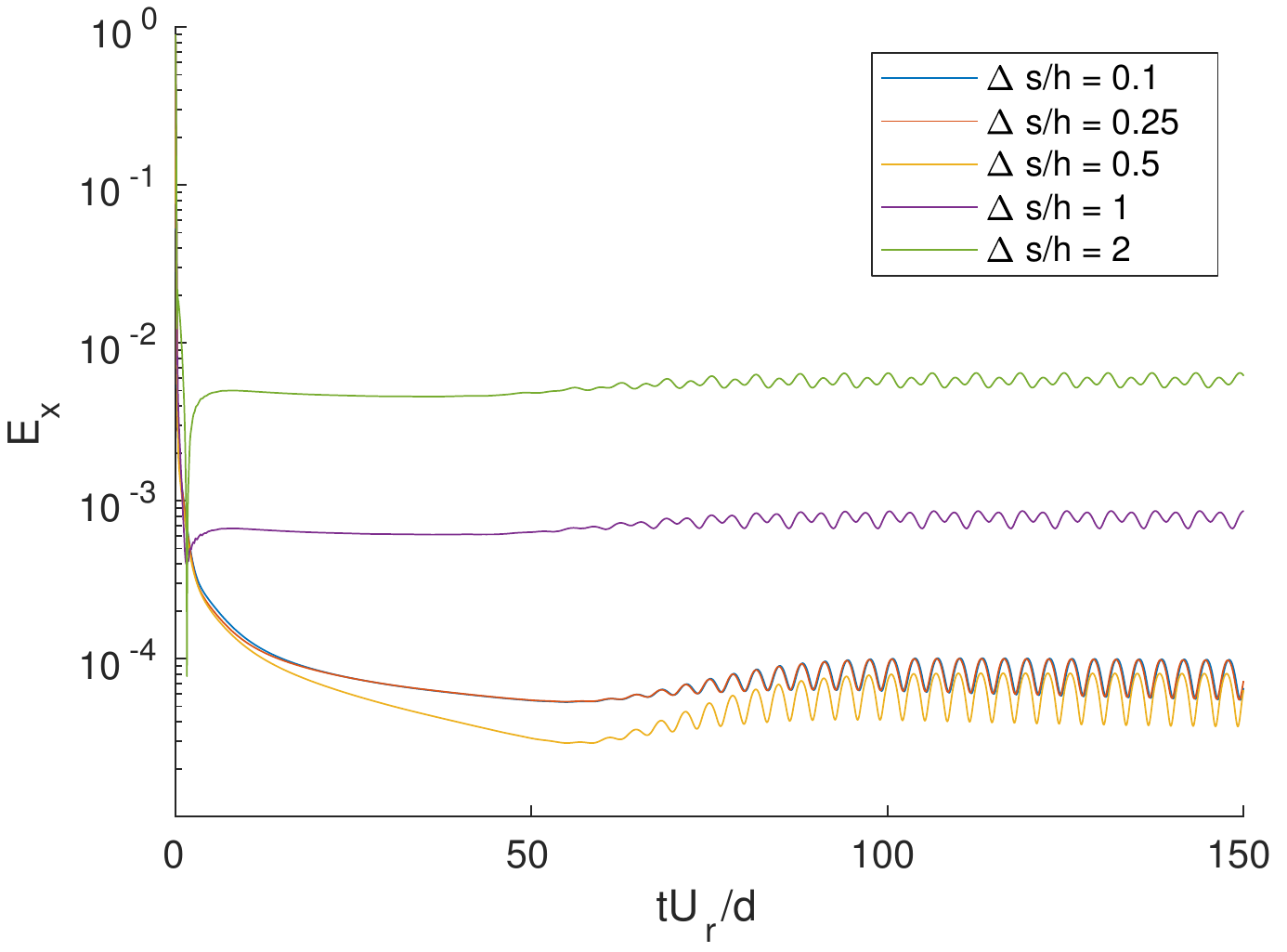}
\put(13,75){\small{BDF1, -$\alpha\Delta t^2$ = 3.9, -$\beta\Delta t$ = 1.9, -$\gamma = 0$}}
      \end{overpic}
 \begin{overpic}[width=0.45\textwidth]{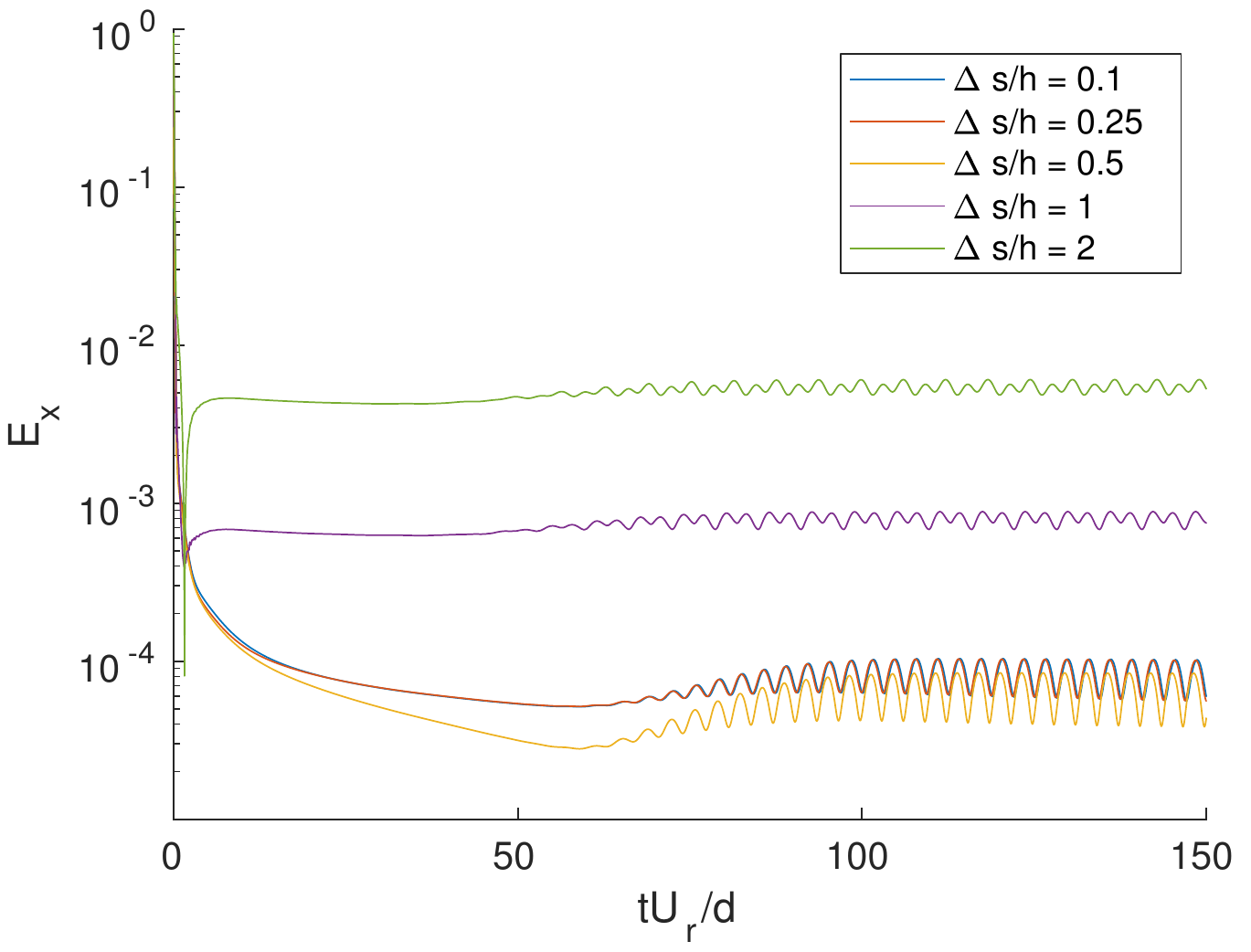}
       \put(16,75){\small{BDF2, -$\alpha\Delta t^2$ = 3.9, -$\beta\Delta t$ = 1.9, -$\gamma = 2$}}
   \end{overpic}\\
\caption{Flow past a stationary cylinder at Re = 100: influence of $\Delta s/h$ on the time evolution $E_x$.}
\label{fig:Err_N_L}
\end{figure}

The first test we consider is the flow past a stationary cylinder at $\mbox{Re} = 100$. Table \ref{tab:mesh} reports details related to the orthogonal Cartesian grids that we use. Concerning meshes with $260\mbox{k}$ and $1000\mbox{k}$ cells, the cylinder has a radius of $0.15$ and its center is located at ($1.85$, $4.0$). On the other hand, for the mesh with $117\mbox{k}$ cells, the cylinder has a radius of $0.5$ and its center is located at ($0$, $0$). 
All the meshes are uniform with the exeption of the mesh with $117\mbox{k}$ cells that however has a uniform local refinement next to the region occupied by the cylinder comparable with the grid spacing of the others meshes. These meshes have been selected to compare our results with computational data reported in \cite{Shin2008}, in which the VBM is combined with a Finite Difference method and the BDF1 time scheme is used. We impose free-stream boundary condition, $\u = \left(1, 0\right)$, at the inflow and far-field boundaries, and an advective boundary condition,
\begin{equation}\label{ref:advective}
\dfrac{\partial \u}{\partial t} + a \dfrac{\partial \u}{\partial \n} = 0,
\end{equation} 
at the outflow, where $a$ is the advection velocity magnitude computed so that the total mass is conserved, and $\n$ is the outwards normal. The partitioned algorithm we use (see Sec. \ref{sec:space_discrete}) requires a boundary condition for the pressure too. We choose $p = 0$ at the outflow and $\partial p/ \partial \n = 0$ on all the other boundaries. We set the density $\rho = 1$ and the viscosity $\mu = 10^{-2}, 3\cdot10^{-3}$ for mesh $117\mbox{k}$, and for meshes $260\mbox{k}$ and $1000\mbox{k}$, respectively. We start the simulations from fluid at rest. The quantities of interest for this benchmark are the drag and lift coefficients given by \cite{Lai2000}, 

\begin{equation}
C_d = -\dfrac{2}{\rho L_r U_r^2} \sum_{i} f_x^i, \quad C_l = -\dfrac{2}{\rho L_r U_r^2} \sum_{i} f_y^i,
\end{equation} 
where $f_x$ and $f_y$ are streamwise and crosswise component of the Euler forcing respectively, $U_r = 1$ is the maximum velocity at far-field boundaries and $L_r = d$ is the diameter of the cylinder, as well as the Strouhal number, $\mbox{St}$, based on the oscillation frequency of the lift force $f$,
\begin{equation}\label{eq:Str}
\mbox{St} = \dfrac{fL_r}{U_r}.
\end{equation}
Table \ref{tab:cyl_staz_1} shows the temporal average of the drag coefficient, $\bar{C}_d$, the amplitude of the lift coefficient oscillation, $C'_l$, and $\mbox{St}$ obtained using the proposed method. Also, computational data provided in \cite{Shin2008}, as well as results related to other studies \cite{Lai2000, Uhlmann2005, Kim2001, Constant2017, Lee2006, Linnick2005, Huang2007_2}, are reported. Notice that when the BDF2 time scheme is used with $\alpha \neq 0$ and $\beta = 0$, ones should necessarily set $\gamma \neq 0$ based on what reported in Sec. \ref{sec:stability_analysis}. 
We observe that \emph{Case 1} and \emph{Case 2}, when the Finite Volume method here proposed is adopted, provide significantly different results with respect to the Finite Difference method used in \cite{Shin2008}. 
Moreover, in \cite{Shin2008}, an improvement of the results is obtained when the domain size is enlarged downstream and on the top (\emph{Case 3}) whilst we do not appreciate significant differences. On the contrary, we are able to obtain better predictions by enlarging the domain size in all the directions (\emph{Case 4}). Results obtained with the two time schemes here considered, BDF1 and BDF2, show a significant difference on the value of $C'_l$ for the mesh $260\mbox{k}$ (\emph{Case 1}). Such a discrepancy decreases when the time step is reduced (\emph{Case 2}), as well as the domain size is enlarged (\emph{Case 3} and \emph{Case 4}). Numerical experiments at maximum time steps based on eq. \eqref{eq:BDF2vsBDF1} for both time schemes are also reported for \emph{Case 4}. Notice that the values of the maximum $CFL$ number obtained by the present simulations are in perfect agreement with ones reported in \cite{Shin2008}. Figure \ref{fig:vel_mag} shows a snapshot of the velocity field and Figure \ref{fig:coeffs} displays the time evolution of the drag and lift coefficients obtained for \emph{Case 4}. We observe that results obtained by using different time schemes show differences in terms of phase. It should be noted that the mesh with $117\mbox{k}$ was not used in \cite{Shin2008} for the case of a stationary cylinder but a very similar mesh was adopted for the simulation of inline/transverse oscillation of a circular cylinder (see Secs. \ref{sec:osc_y_cyl} and \ref{sec:osc_x_cyl}).
For further comparison, we evaluate also the error for the virtual boundary velocity in the streamwise direction, $E_x$, defined as

\begin{align}\label{eq:error_x}
E_x = \left(\dfrac{1}{N_L} \sum_{k = 1}^{N_L} \left((u_{ib})_k - (u_{b})_k\right)^2\right)^{1/2}.
\end{align}
We are going to investigate the influence of the feedback forcing gains, $\alpha$, $\beta$ and $\gamma$, as well as the number of Lagrangian points, $N_L$, on the error, for the mesh $117\mbox{k}$. We select the same time step, $\Delta t$ = $0.01$, of \cite{Shin2008}, and the values of $-\alpha \Delta t^2$ and $-\beta \Delta t$ include, among others, those ones used in \cite{Shin2008}. Moreover, concerning the tests carried out by varying the forcing parameters, we set $\Delta s/h = 1$ as done in \cite{Shin2008}. However, we observe that in \cite{Shin2008} such an analysis was carried out for $\mbox{Re} = 185$ by using a mesh with about the half of the cells and a 3-point regularized delta function. Figure \ref{fig:Err_x_alpha} shows the evolution over time of $E_x$ varying $-\alpha \Delta t^2$ for $-\beta \Delta t = 1.5$ and $-\gamma = 0$. We observe that the error converges to a smaller value for the larger value of  $-\alpha \Delta t^2$. This result is in agreement with that one in \cite{Shin2008}. 
We note that the difference between $-\alpha \Delta t^2 = 1$ and $-\alpha \Delta t^2 = 3.5$ is very little. This suggests that the convergence is reached for low values of $-\alpha \Delta t^2$, far from the stability limits (see Sec. \ref{sec:stability_analysis}). Figure \ref{fig:Err_x_beta} shows the evolution over time of $E_x$ varying $-\beta \Delta t$ for $-\alpha \Delta t^2 = 0.1$ and $-\gamma = 0$. We observe that the transient decay of the error is greater for the larger value of $-\beta \Delta t$. On the other hand, all the tested cases show the same level of convergence of the error, suggesting that $-\alpha \Delta t^2$ is a more critical parameter. These results are in agreement with \cite{Shin2008}. Nevertheless, in our calculations, $-\beta \Delta t$ does not reduce the convergence time, as showed in \cite{Shin2008}. Moreover, errors are not in phase as in \cite{Shin2008}.
Figure \ref{fig:Err_x_gamma} shows the evolution over time of $E_x$ varying $-\gamma$ for $-\alpha \Delta t^2 = 0.1$ and $-\beta \Delta t = 1.5$. We observe a very low sensitivity with respect to $-\gamma$. Figure \ref{fig:Err_N_L} shows the evolution over time of $E_x$ varying $\Delta s/h$. We observe a non monotonic convergence of the error when $\Delta s/h < 1$. 
Finally, we observe that the sensitivity of the error with respect to the time scheme is very low.

For the sake of completeness, Tables \ref{tab:cyl_staz_alpha}, \ref{tab:cyl_staz_beta}, \ref{tab:cyl_staz_gamma} and \ref{tab:cyl_staz_ds_h} show the influence of the feedback forcing gains, $\alpha$, $\beta$ and $\gamma$, and of the number of Lagrangian points, $N_L$, on the values of $\bar{C}_d$, $C'_l$ and $\mbox{St}$. We highlight that, as in \cite{Shin2008}, the number of Lagrangian points has a negligible effect on the results. Moreover, it is confirmed that $\alpha$ is the forcing parameter that mostly affect the performances of the method. Results related to the BDF2 time scheme, not reported for brevity, does not show substantial difference with respect to this scenario. 

We notice that, despite the fact that a Finite Volume approximation is used in \cite{Lee2006}, a detailed comparison with our results is complicated because computational features, such as time step, number of Lagrangian points, errors, are missing. 

In conclusion, we learned that for this benchmark $-\alpha\Delta t^2$ is the most critical parameter to be properly tuned in order to optimize the accuracy and efficiency of the computation. 

\subsection{Oscillatory flow past a stationary cylinder}\label{sec:John}
In order to provide a further test of the capabilities of the proposed method to simulate flow fields involving bodies at rest, we additionally consider a periodic oscillation of the flow past a stationary cylinder \cite{John2004, turek1996, Girfoglio2019}. The computational domain is a $2.2 \times 0.41$ rectangular channel with a cylinder of radius $0.05$ centered at ($0.2$, $0.2$), where the bottom left corner of the channel as the origin of the axes. We impose a no slip boundary condition on the upper and lower wall. At the inflow and the outflow we prescribe the following velocity profile:
\begin{align}\label{eq:cyl_bc}
\u(0,y,t) = \left(\dfrac{6}{0.41^2} \sin\left(\pi t/8 \right) y \left(0.41 - y \right), 0\right), \quad y \in [0, 2.2], \quad t \in (0, 8].
\end{align}

\begin{figure}
\centering
 \begin{overpic}[width=0.45\textwidth]{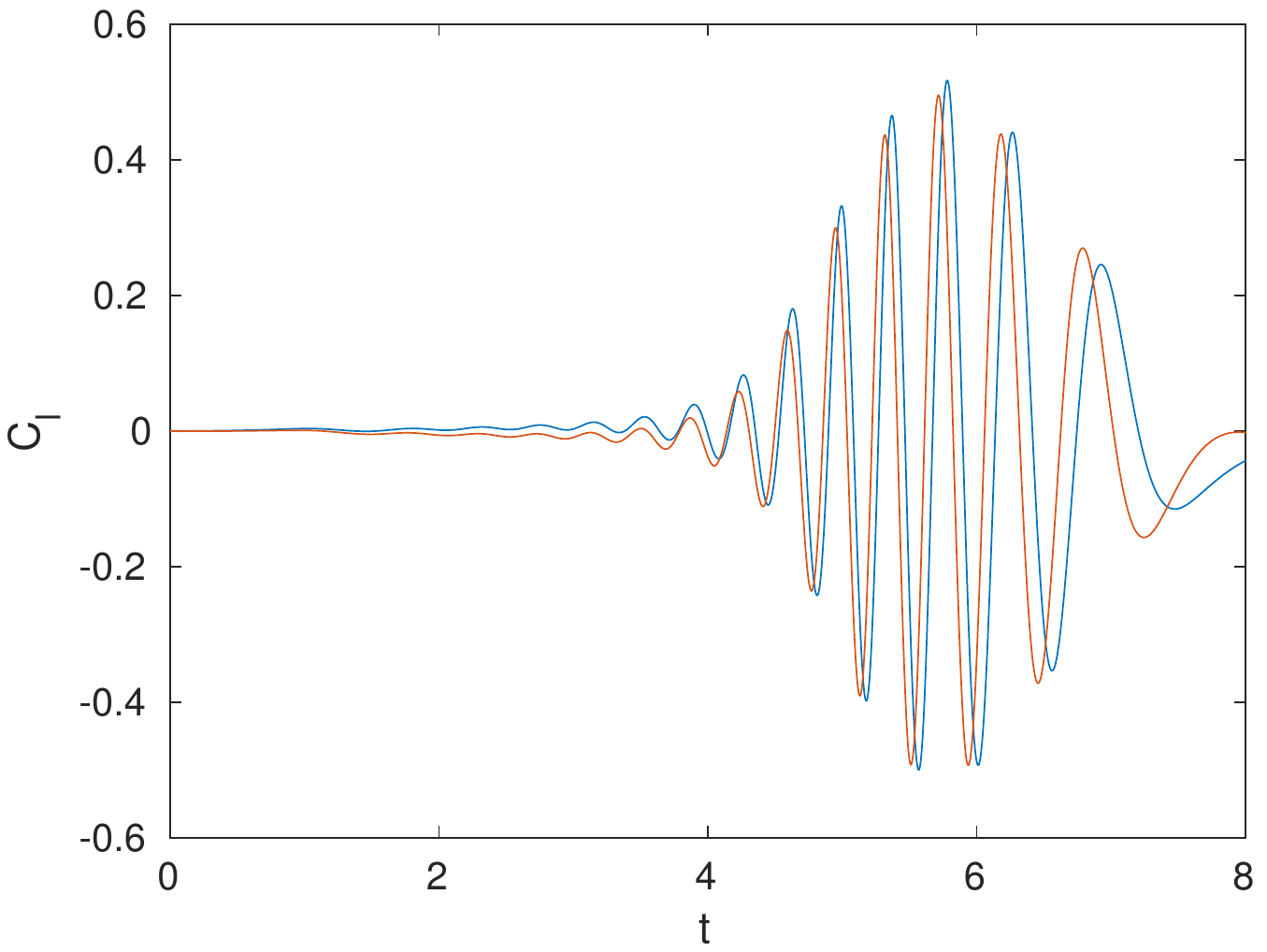}
      \end{overpic}
 \begin{overpic}[width=0.45\textwidth]{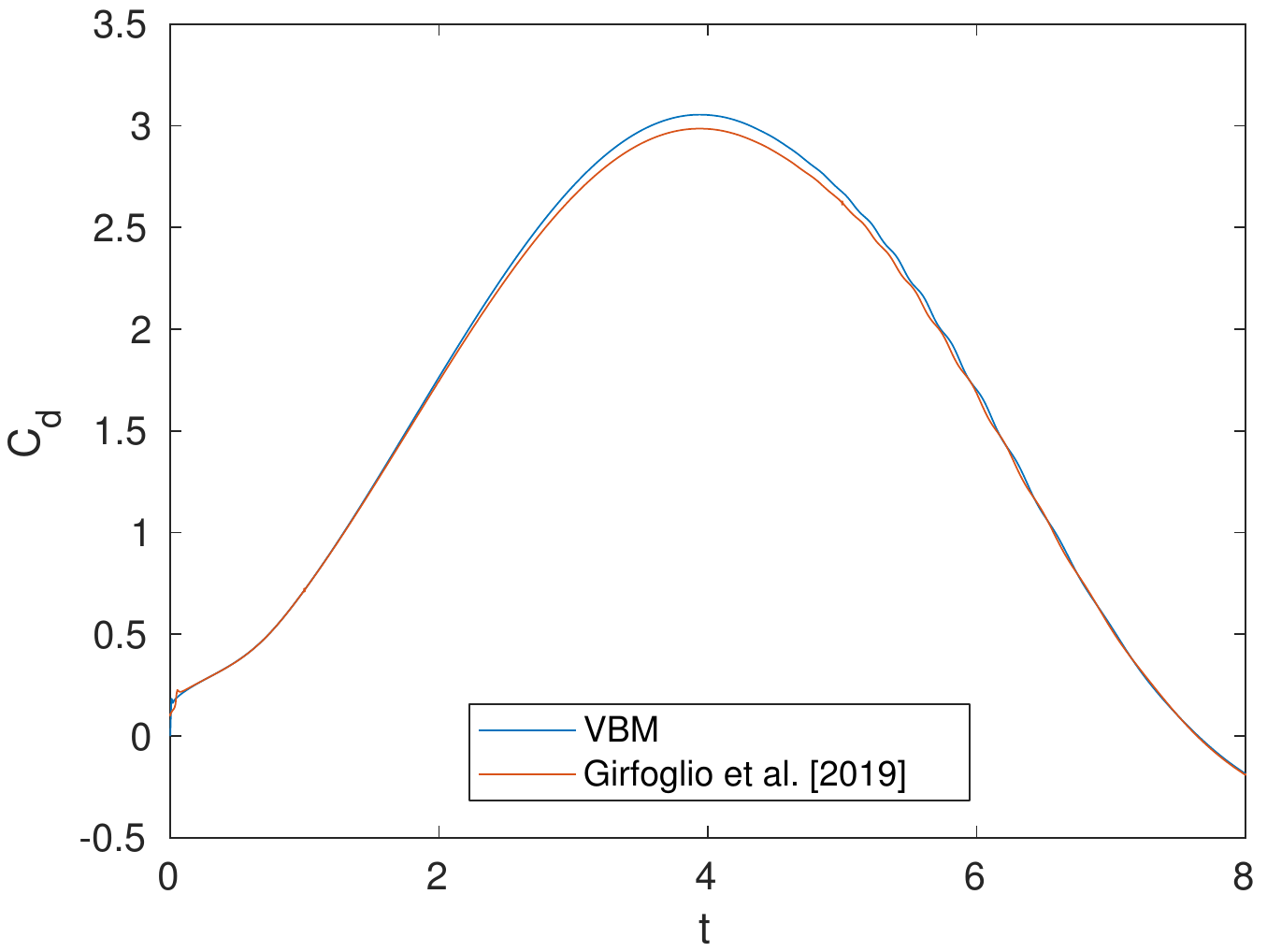}
   \end{overpic}\\
\caption{Oscillatory flow past a stationary cylinder: time evolution of lift and drag coefficients against the results reported in \cite{Girfoglio2019}.}
\label{fig:John2004}
\end{figure}

Concerning the pressure, we impose $\partial p/\partial \n = 0$ on all the boundaries. We set density $\rho = 1$ and viscosity $\mu = 10^{-3}$. We start the simulations from fluid at rest. Note that the Reynolds number is time dependent, with $0 \leq \mbox{Re} \leq 100$ \cite{John2004,turek1996,Girfoglio2019}. 
Such benchmark requires a roughly uniform hexaedral/prismatic mesh of about $200\mbox{k}$ cells for a DNS with a standard Finite Volume method \cite{Girfoglio2019}. Here we use a uniform orthogonal Cartesian mesh of $200\mbox{k}$ cells. Based on the results in Sec. \ref{sec:staz_cyl}, we set $-\alpha \Delta t^2 = 1$ and $-\beta \Delta t = -\gamma = 0$, $\Delta s/h = 1$, and we limit to consider only the BDF1 time scheme. Figure \ref{fig:John2004} shows the time evolution of the lift and drag coefficients that are in very good agreement with those ones reported in \cite{Girfoglio2019}.

\subsection{Transverse oscillation of a circular cylinder in a free-stream}\label{sec:osc_y_cyl}
\begin{figure}
\centering
 \begin{overpic}[width=0.45\textwidth]{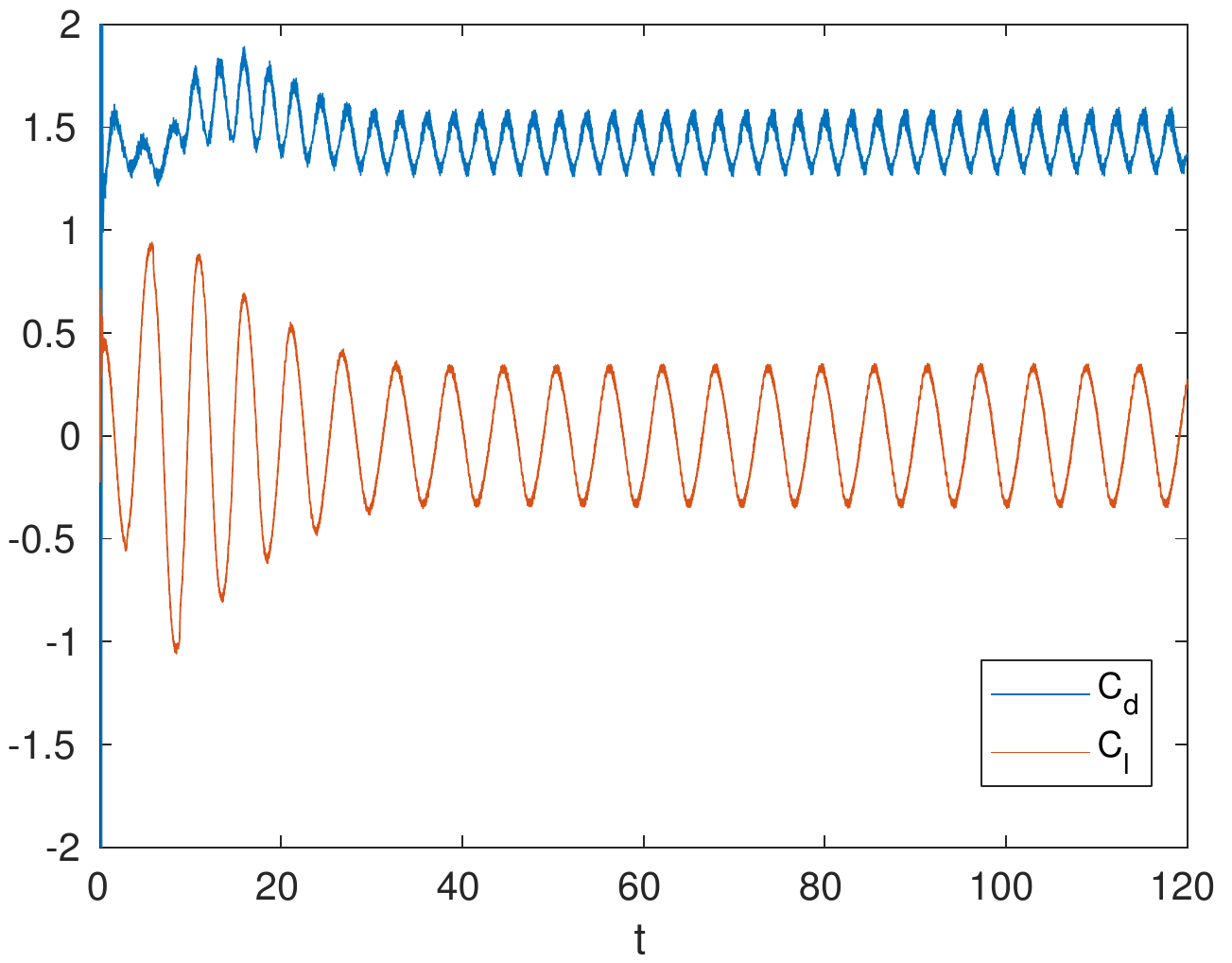}
       \put(48,80){\small{a)}}
      \end{overpic}
 \begin{overpic}[width=0.45\textwidth]{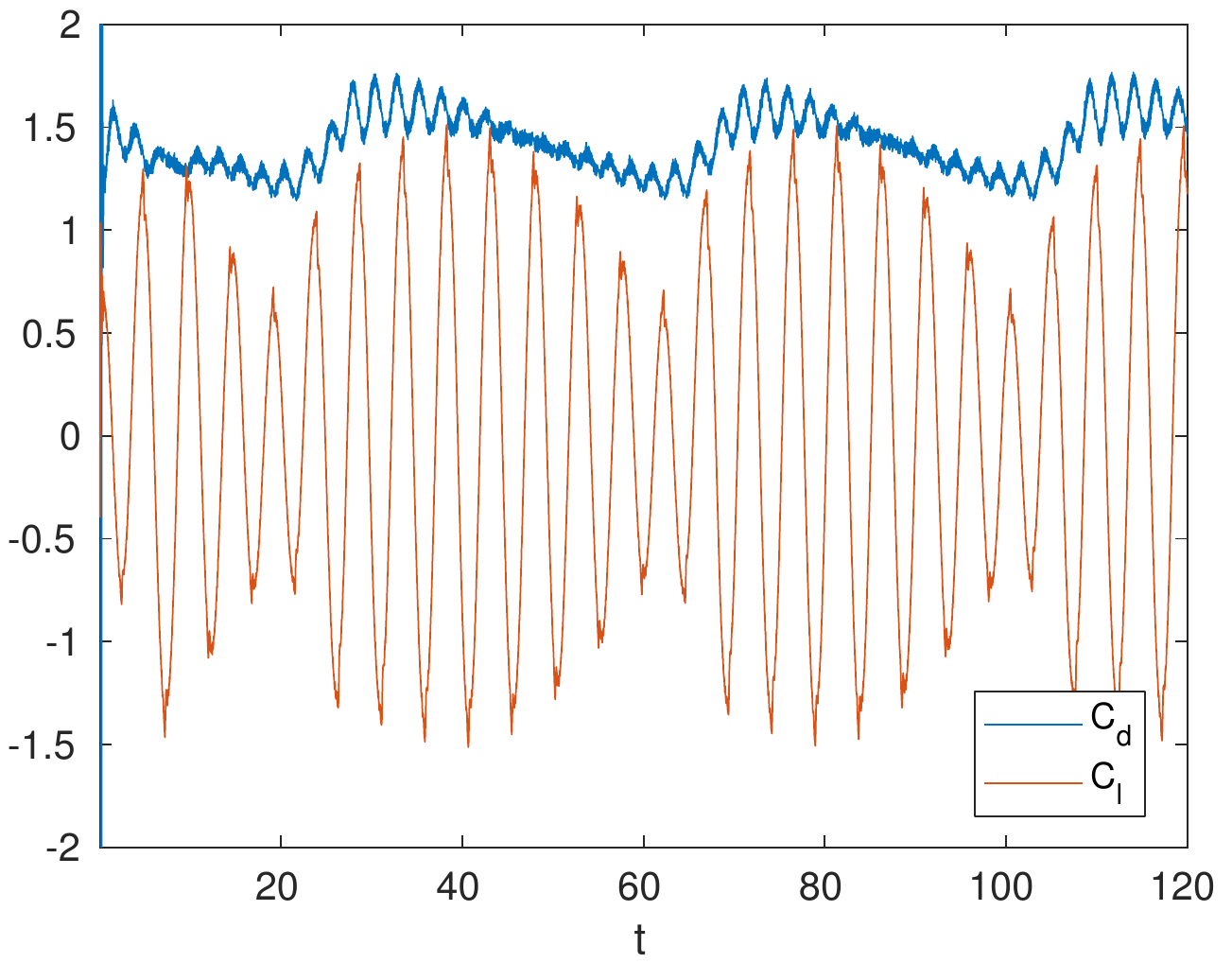}
       \put(48,80){\small{b)}}
   \end{overpic}\\
\caption{Transverse oscillation of a circular cylinder in a free-stream: a) $f_e/f_0 < 1$, b) $f_e/f_0 > 1$.}
\label{fig:mov_y_1}
\end{figure}

Next, we consider a benchmark involving a body that rigidly moves. We investigate a periodic transerve oscillation of a circular cylinder in a free-stream \cite{Guilmineau2002}. The time-periodic motion of the center of the cylinder is given by 

\begin{equation}\label{eq:y_osc}
\left(x_c(t),y_c(t)\right) = \left(0, A_{m} \cos \left(2\pi f_e t\right)\right),
\end{equation}
where $A_m$ and $f_e$ are the amplitude and the frequency of the oscillation, respectively. The Reynolds number, based on the free-stream velocity, is 185, and $A_m/d = 0.2$ \cite{Guilmineau2002}. We impose free-stream boundary condition, $\u = (1, 0)$, at the inflow and far-field boundaries, and advective boundary condition \eqref{ref:advective} at the outflow. For comparison with \cite{Shin2008}, we consider the mesh $117\mbox{k}$ and perform a first simulation by setting $-\alpha\Delta t^2 = 3.9$, $-\beta\Delta t = 1.9$, $-\gamma = 0$ and $\Delta s/h = 1$. We investigate the two possible cases, $f_e/f_0 < 1$ and $f_e/f_0 > 1$, where $f_0 = 0.19$ is the natural shedding frequency for a stationary cylinder at $\mbox{Re} = 185$ \cite{Guilmineau2002}. The computational time step, $\Delta t = T/720$, is based on the period of the oscillation $T = 1/f_e$, and the BDF1 time scheme is employed. Finally, we set density $\rho = 1$ and viscosity $\mu = 5.4e-3$. Figure \ref{fig:mov_y_1} shows time evolution of the drag and lift coefficients for $f_e/f_0 < 0.9$ (a) and $f_e/f_0 > 1$ (b). We can note that the coefficients shows a regular trend once vortex shedding is established. Moreover, we observe that for $f_e/f_0 > 1$, both the drag and lift coefficients exhibit beat phenomena. These results are in good agreement with those ones reported in \cite{Guilmineau2002, Shin2008}. 

\begin{figure}
\centering
 \begin{overpic}[width=0.45\textwidth]{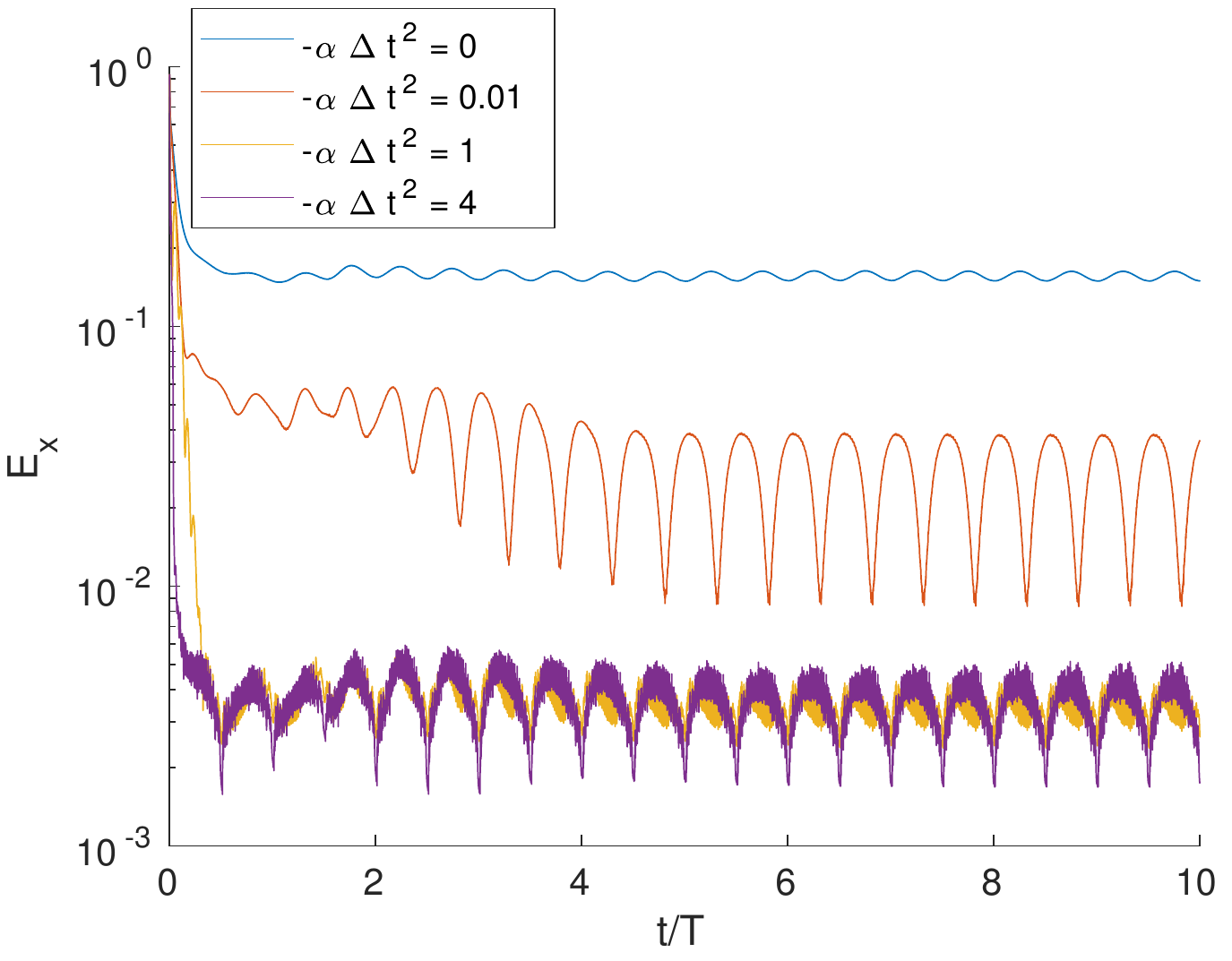}
       \put(32,80){\small{BDF1, -$\beta\Delta t$ = 1, $\gamma$ = 0}}
      \end{overpic}
 \begin{overpic}[width=0.45\textwidth]{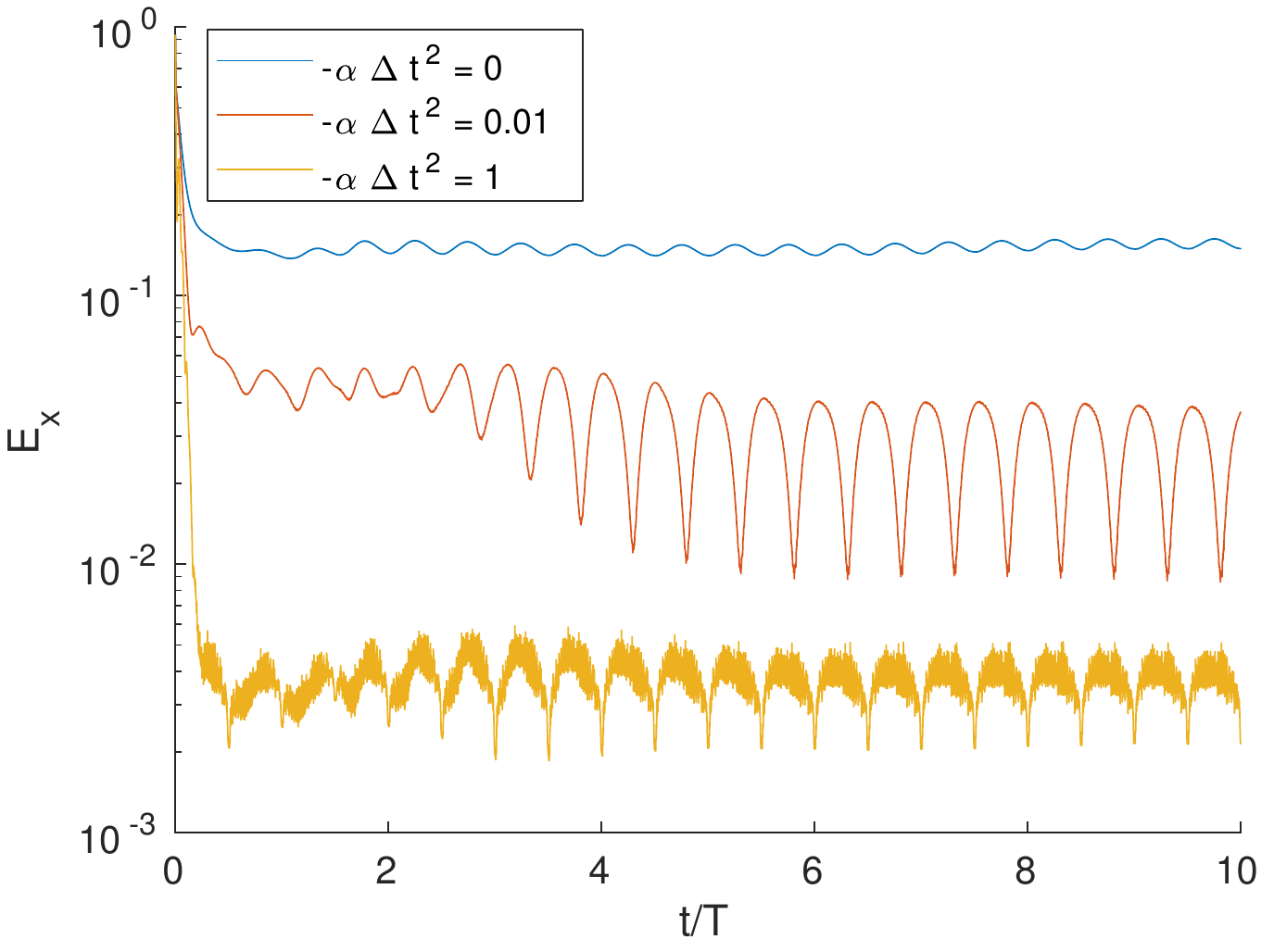}
       \put(32,80){\small{BDF2, -$\beta\Delta t$ = 1, $\gamma$ = 0}}
   \end{overpic}\\
\caption{Transverse oscillation of a circular cylinder in a free-stream: influence of $-\alpha \Delta t^2$ on the time evolution of $E_x$.}
\label{fig:Err_x_alpha_mov_y}
\end{figure}

\begin{figure}
\vspace{3cm}
\centering
 \begin{overpic}[width=0.45\textwidth]{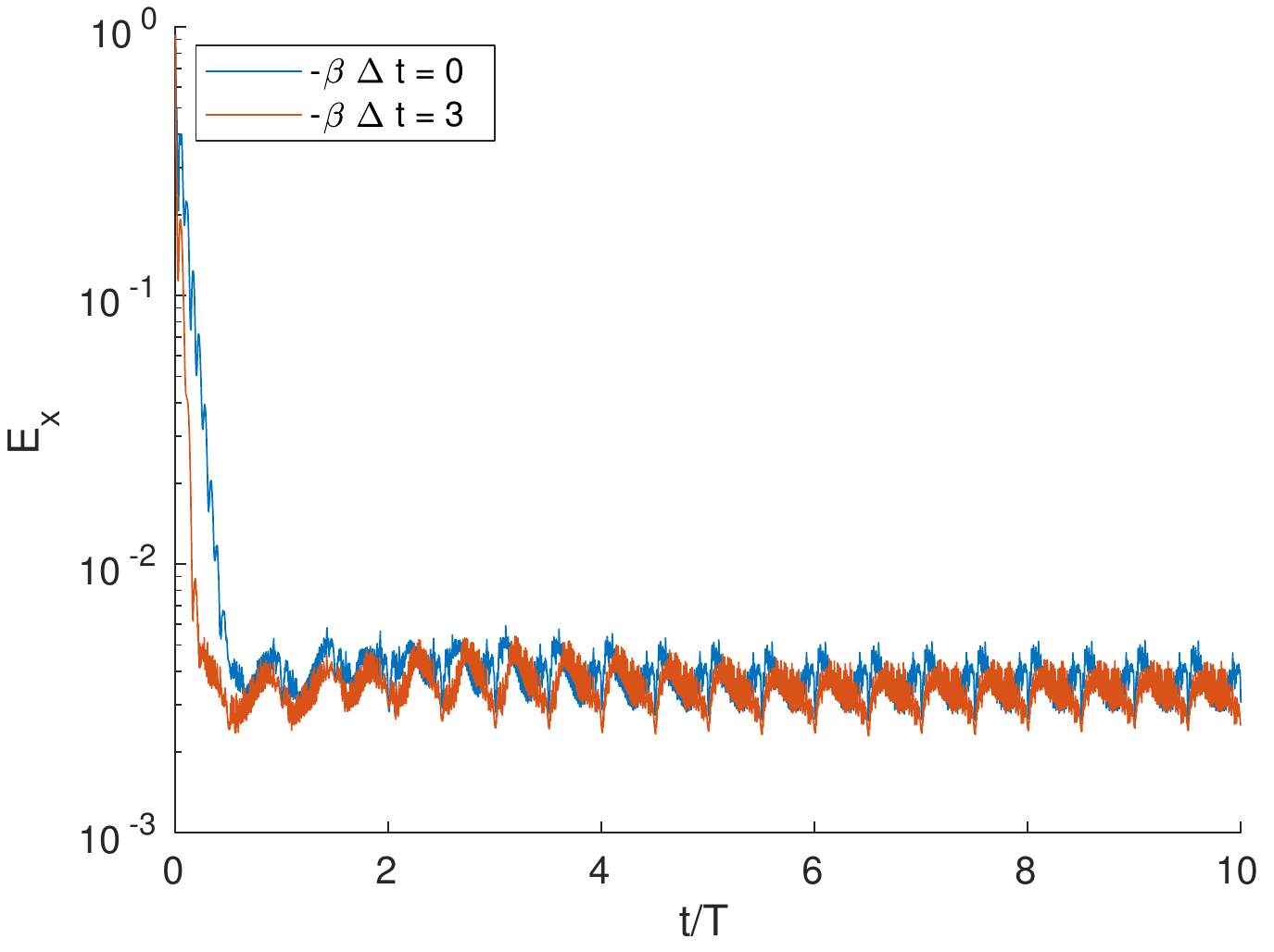}
       \put(32,75){\small{BDF1, -$\alpha\Delta t^2$ = 0.4, $\gamma$ = 0}}
      \end{overpic}
 \begin{overpic}[width=0.45\textwidth]{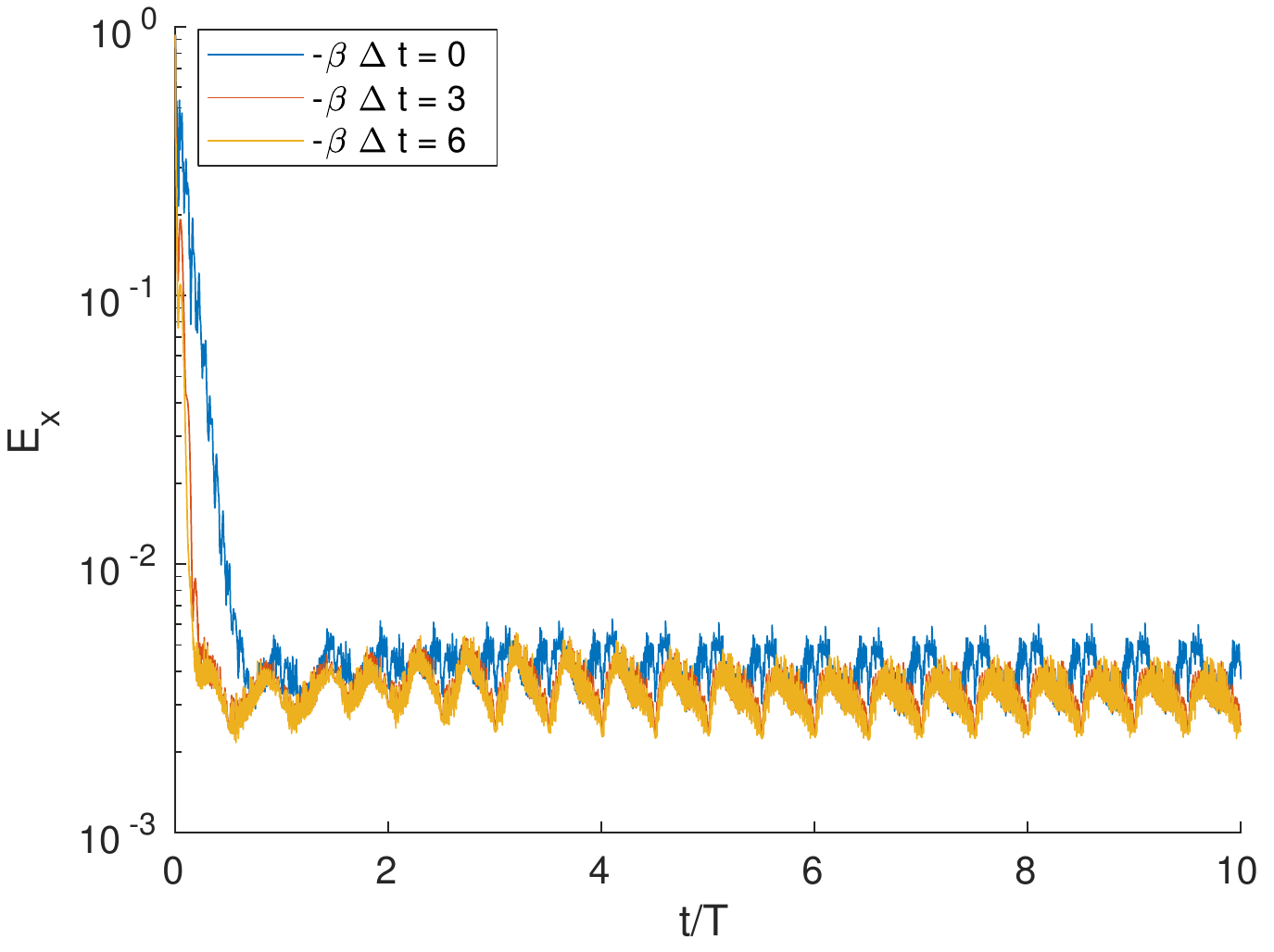}
       \put(32,75){\small{BDF2, -$\alpha\Delta t^2$ = 0.4, $\gamma$ = 0}}
   \end{overpic}\\
\caption{Transverse oscillation of a circular cylinder in a free-stream: influence of $-\beta \Delta t$ on the time evolution of $E_x$.}
\label{fig:Err_x_beta_mov_y}
\end{figure}

\begin{figure}
\centering
 \begin{overpic}[width=0.45\textwidth]{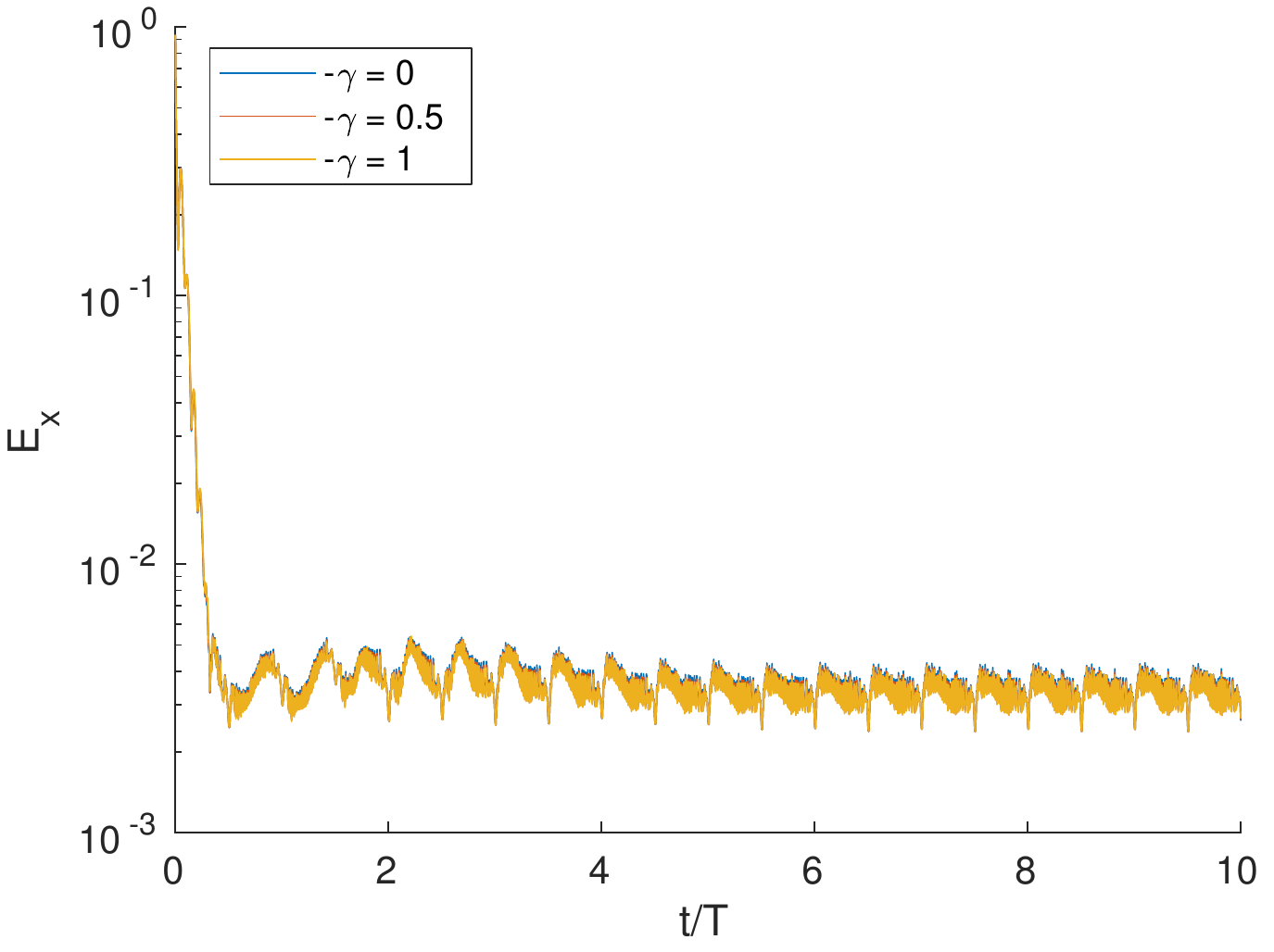}
       \put(32,73){\small{BDF1, -$\alpha\Delta t^2$ = 0.4, -$\beta\Delta t$ = 1}}
      \end{overpic}
 \begin{overpic}[width=0.45\textwidth]{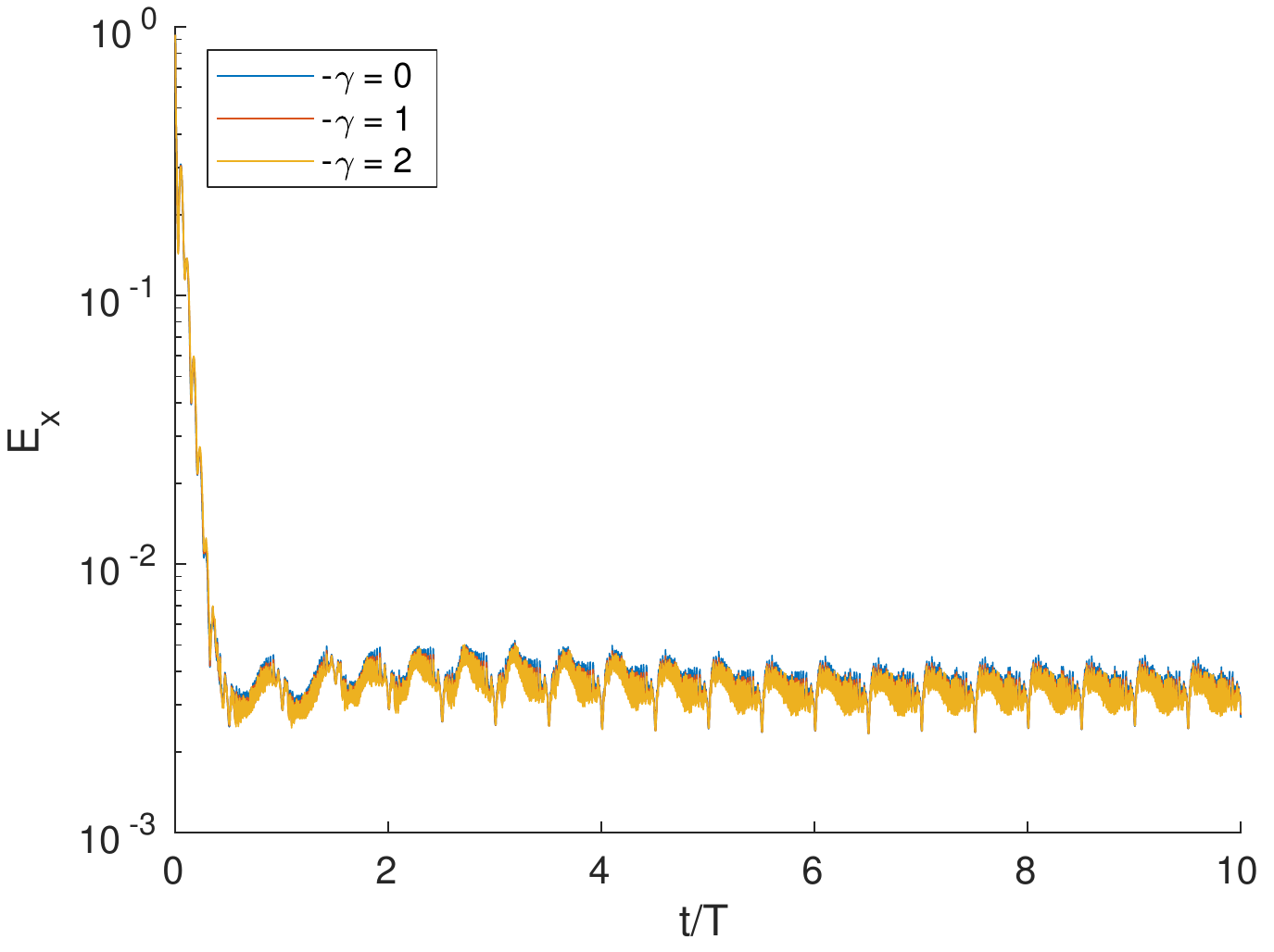}
       \put(32,73){\small{BDF2, -$\alpha\Delta t^2$ = 0.4, -$\beta\Delta t$ = 1}}
   \end{overpic}\\
\caption{Transverse oscillation of a circular cylinder in a free-stream: influence of $-\gamma$ on the time evolution of $E_x$.}
\label{fig:Err_x_gamma_mov_y}
\end{figure}

\begin{figure}
\centering
 \begin{overpic}[width=0.45\textwidth]{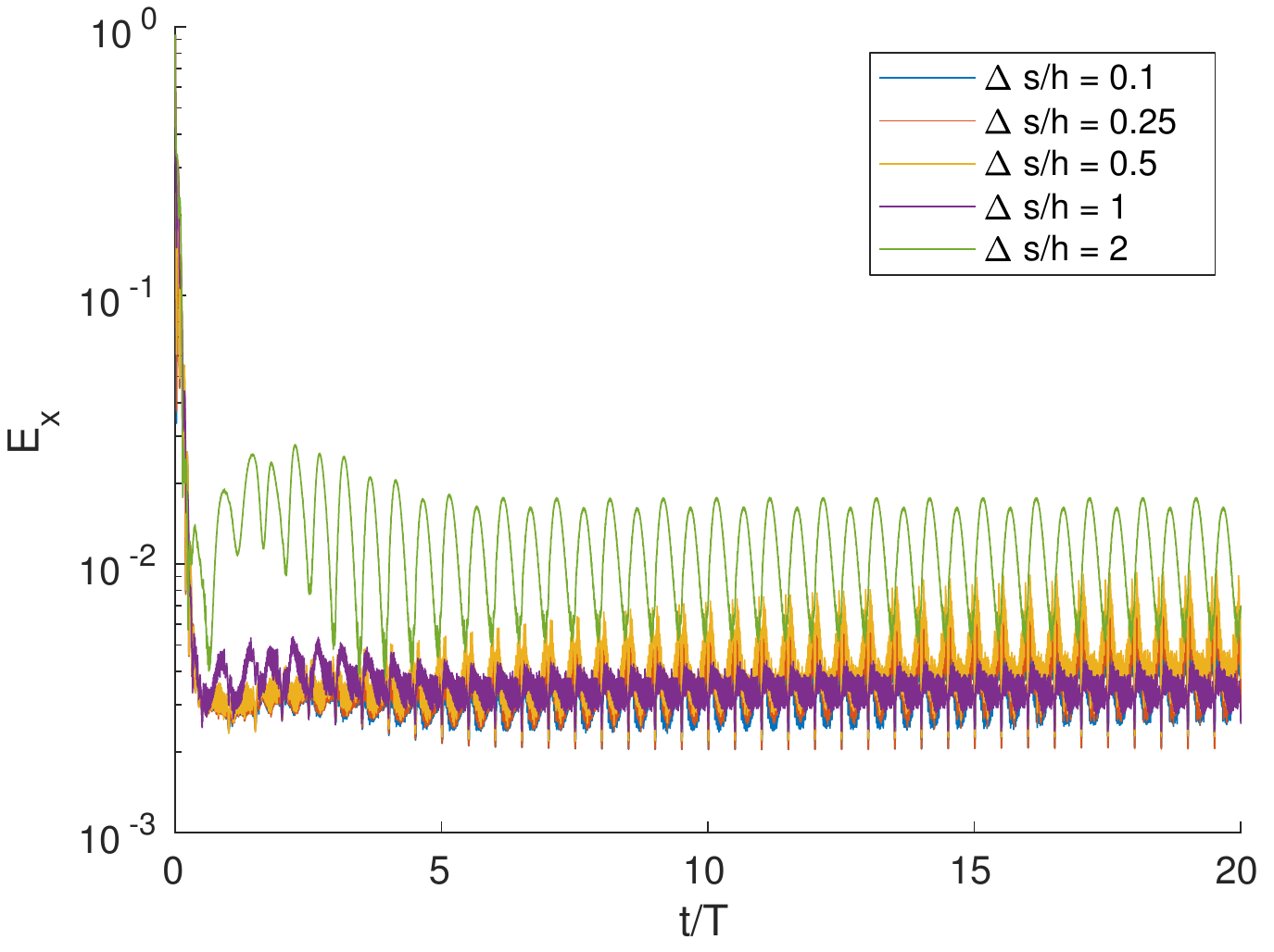}
       \put(16,75){\small{BDF1, -$\alpha\Delta t^2$ = 0.4, -$\beta\Delta t$ = 1, -$\gamma = 1$}}
      \end{overpic}
 \begin{overpic}[width=0.45\textwidth]{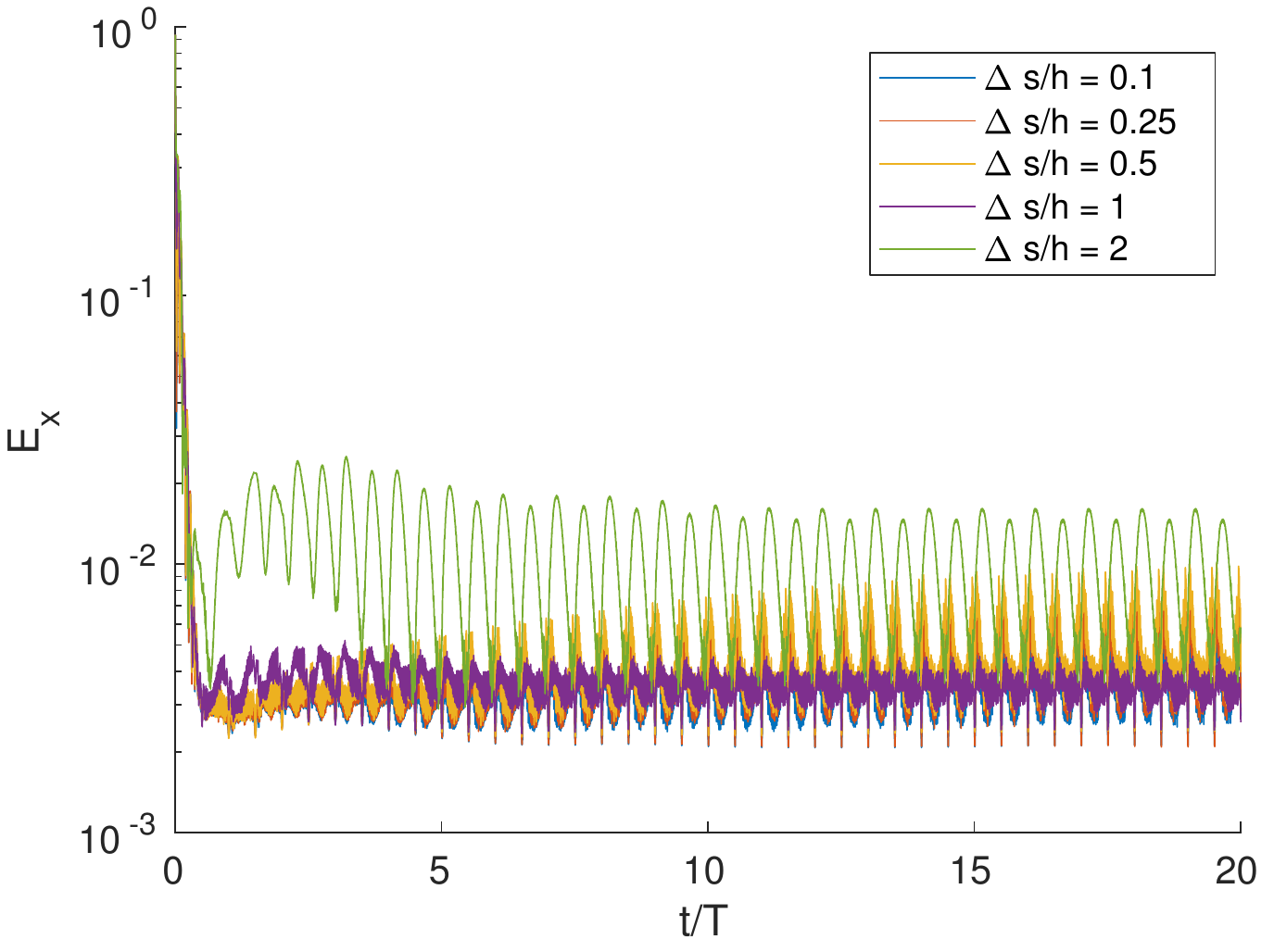}
       \put(19,75){\small{BDF2, -$\alpha\Delta t^2$ = 0.4, -$\beta\Delta t$ = 1, -$\gamma = 1$}}
   \end{overpic}\\
\caption{Transverse oscillation of a circular cylinder in a free-stream: influence of $\Delta s/h$ on the time evolution of $E_x$.}
\label{fig:Err_x_N_L_mov_y}
\end{figure}

Concerning the investigation of the error, we focus on $f_e/f_0 < 0.9$. We observe that in \cite{Shin2008} the error analysis was carried out for a mesh having twice the number of the cells, i.e. the half of the grid size, of the mesh $260\mbox{k}$. Figure \ref{fig:Err_x_alpha_mov_y} shows the evolution over time of $E_x$ varying $-\alpha \Delta t^2$ for $-\beta\Delta t$ = 1 and $-\gamma$ = 0. We oberve that the error converges to a smaller value for the larger value of $-\alpha \Delta t^2$ as observed in Sec. \ref{sec:staz_cyl} for problems involving stationary bodies. This result is in agreement with that one in \cite{Shin2008}. We note there is very little difference between $-\alpha \Delta t^2 = 1$ and $-\alpha \Delta t^2 = 4$. This suggests that also for moving problems the convergence is reached for low values of $-\alpha \Delta t^2$, far from the stability limits (See Sec. \ref{sec:stability_analysis}). High frequency spurious oscillations affect the error at increasing of $-\alpha \Delta t^2$. Figure \ref{fig:Err_x_beta_mov_y} shows the evolution over time of $E_x$ for different values of $-\beta \Delta t$ for $-\alpha\Delta t^2$ = 0.4 and $-\gamma$ = 0. We observe that the transient decay of the error is greater for the larger value of $-\beta \Delta t$ but unlike stationary problems the error also converges to a slightly smaller value for larger $-\beta \Delta t$. On the other hand, such effect seems to be stronger in \cite{Shin2008} where by moving from $-\beta \Delta t = 0$ to $-\beta \Delta t = 1$ one obtains a greater reduction of the error with respect to that one obtained in the present study by moving from $-\beta \Delta t = 0$ to $-\beta \Delta t = 3$. Figure \ref{fig:Err_x_gamma_mov_y} shows the evolution over time of $E_x$ varying -$\gamma$ for $-\alpha\Delta t^2$ = 0.1 and $-\beta\Delta t$ = 1. Just like for the stationary problems, we observe a very low sensitivity to -$\gamma$.
Figure \ref{fig:Err_x_N_L_mov_y} shows the evolution over time of $E_x$ varying $\Delta s/h$ for $-\alpha\Delta t$ = 0.4, $-\beta \Delta t$ = 1, and $-\gamma = 0$. We do not observe convergence for the values of $\Delta s/h$ here considered but however, like stationary problems, it is evident that the trend of the error is not monotonic. 
Finally, like stationary problems, we observe with respect to the time scheme.

\begin{figure}
\centering
 \begin{overpic}[width=0.45\textwidth]{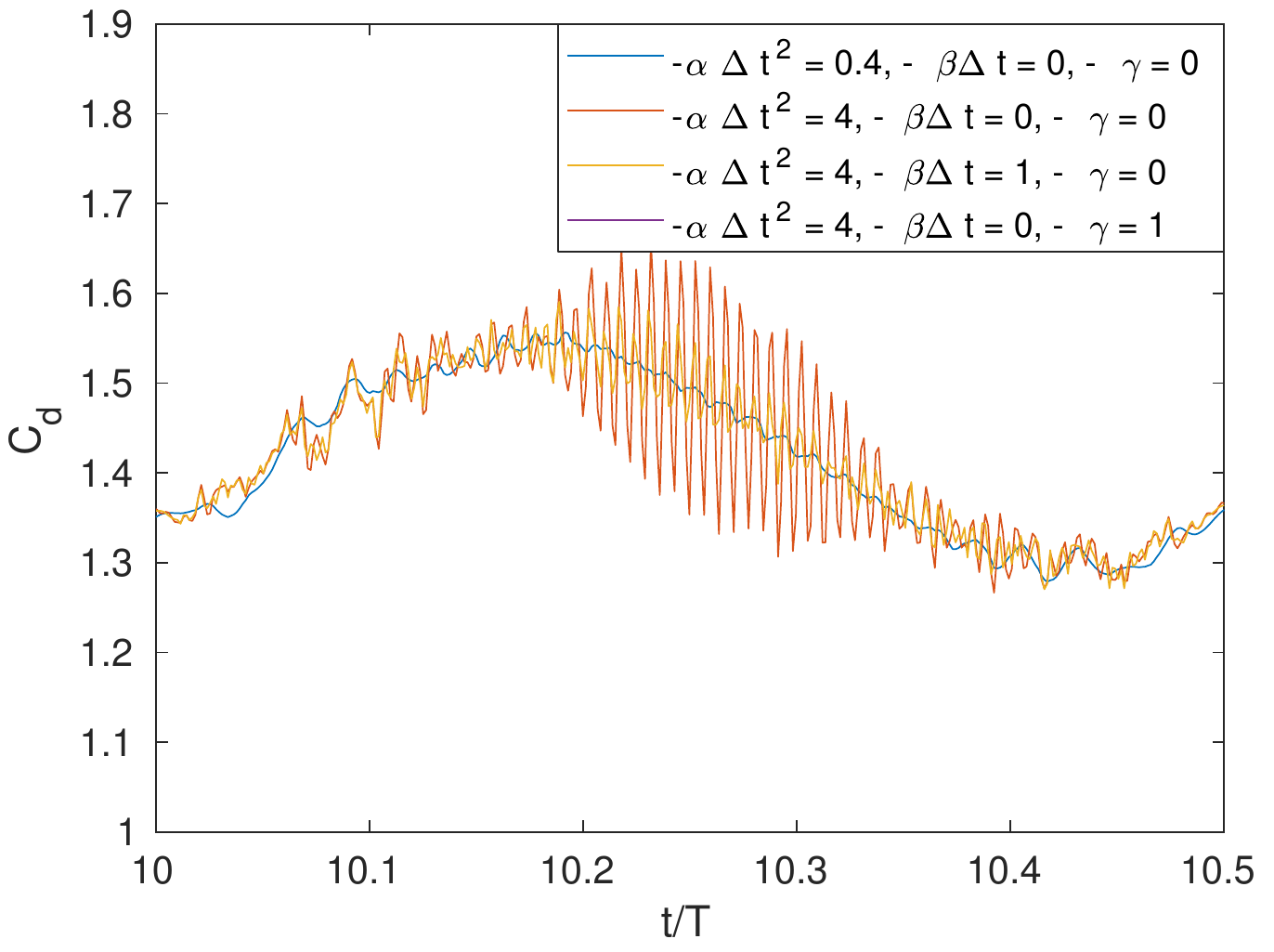}
       \put(40,75){\small{BDF1}}
      \end{overpic}
 \begin{overpic}[width=0.45\textwidth]{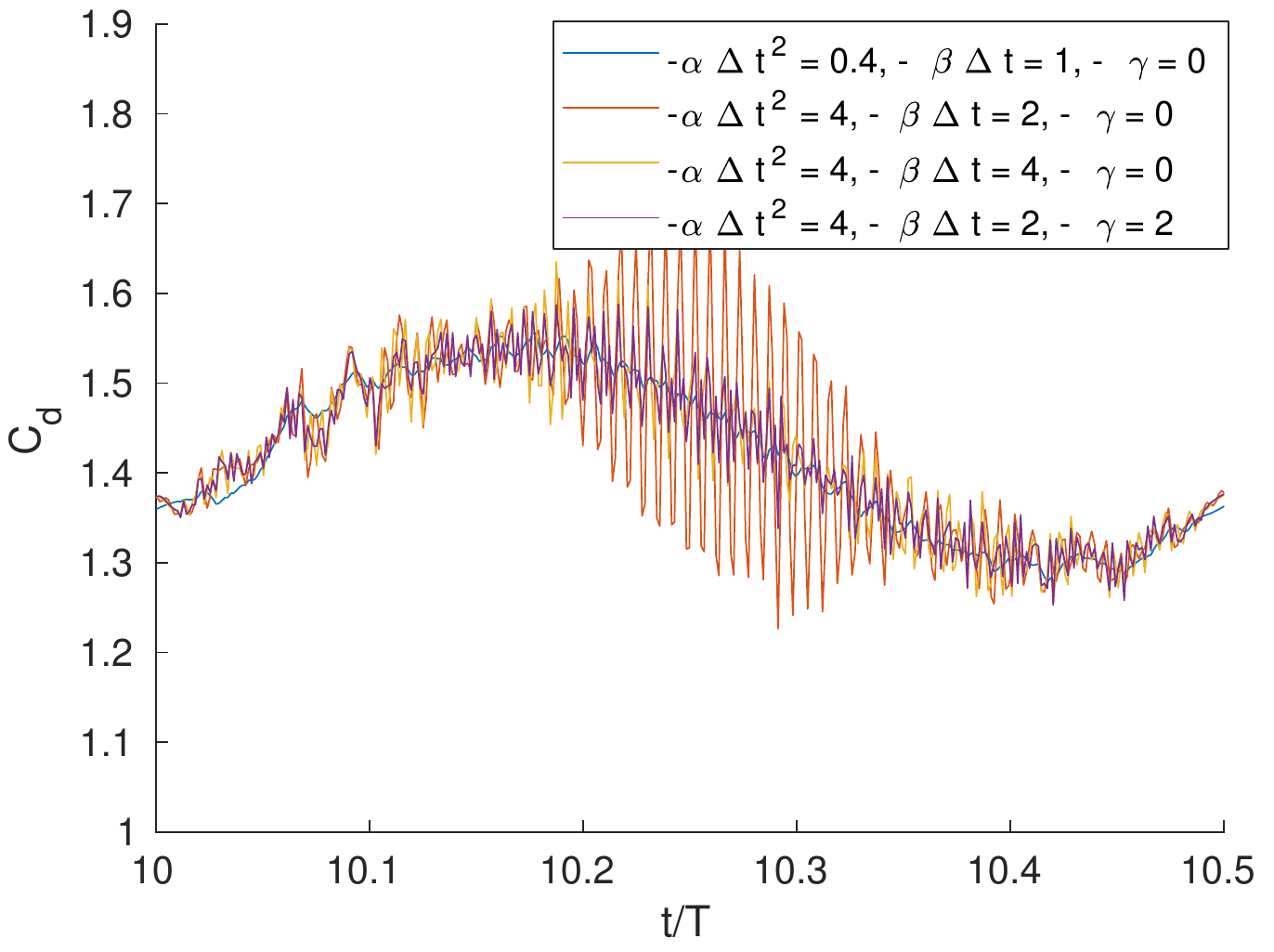}
       \put(40,75){\small{BDF2}}
   \end{overpic}\\
\caption{Transverse oscillation of a circular cylinder in a free-stream: close-up of the time evolution of the drag coefficient for different forcing gains.}
\label{fig:non_grow_osc_y}
\end{figure}

Figure \ref{fig:non_grow_osc_y} displays a close-up of the time history of the drag coefficient for different forcing gains. We observe that $-\beta \Delta t$ helps to reduce spurious oscillations that affect the solution at large $-\alpha \Delta t$. This result is in agreement with that one in \cite{Shin2008}. Moreover, we see that the same effect can be obtained by using $-\gamma$. We note that when the BDF2 time scheme is used a larger value of $-\beta \Delta t$ is necessary to break down the numerical oscillations with respect to the BDF1 time scheme.

In conclusion, we learned that, like stationary problems, even for rigidly moving bodies, $-\alpha\Delta t^2$ is the most critical parameter. Nevertheless, both $-\beta\Delta t$ and $-\gamma$ help to optimize the accuracy and efficiency of the computation.



\subsection{Inline oscillation of a circular cylinder in a quiescent environment}\label{sec:osc_x_cyl}
\begin{figure}
\centering
\includegraphics[width=0.5\textwidth]{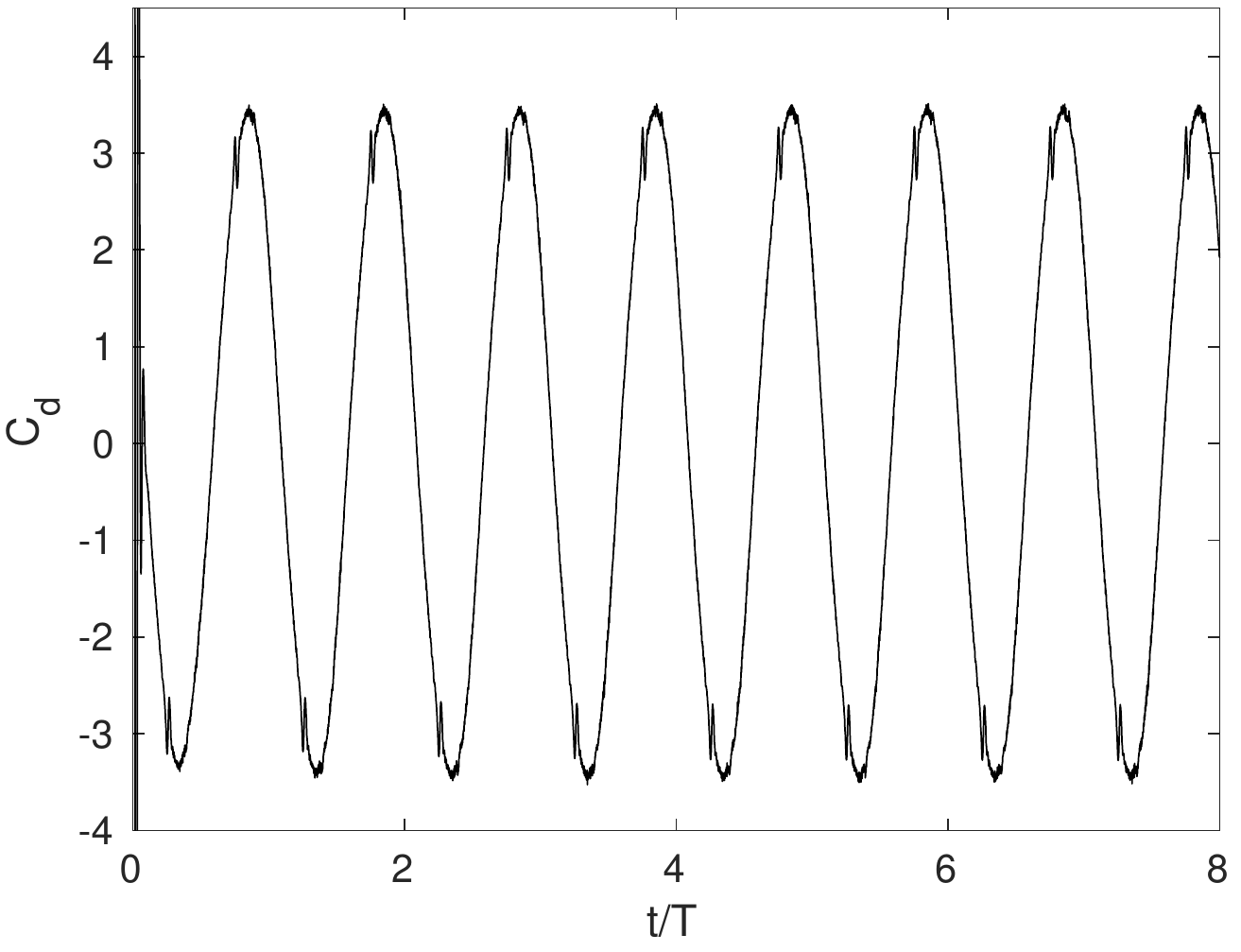}
\caption{Inline oscillation of a circular cylinder in a quiescent environment: time evolution of the drag coefficient at $\mbox{Re} = 100$ and $KC = 5$.}
\label{fig:Cd_osc_x}
\end{figure}

Finally, we consider a periodic inline oscillation of a circular cylinder in a fluid at rest \cite{Dutsch1998}. The time-periodic motion of the center of the cylinder is given by 

\begin{equation}\label{eq:x_osc}
\left(x_c(t), y_c(t)\right) = \left(-A_{m} \sin \left(2\pi f_e t\right),0\right),
\end{equation}
where $A_m$ and $f_e$ are the amplitude and the frequency of the oscillation, respectively. The dimensionless parameters that characterize the dynamics are the Reynolds number, $\mbox{Re}$, and the Keulegan-Carpenter number, $KC$, defined as,

\begin{equation}\label{eq:KC}
KC = \dfrac{U_r}{f_e L_r},
\end{equation}
where $U_r = 2\pi f_e A_m$ and $L_r = d$. We set $\mbox{Re} = 100$ and $KC = 5$, according to \cite{Shin2008,Dutsch1998}. We prescribe do-nothing boundary conditions at all far-fields boundaries. We set $\rho = 1$, $\mu = 10^{-2}$, $A_m = 0.796$ and $f_e = 0.2$. We adopt a computational set-up similar to that one of \cite{Shin2008} to compare the results. We use the mesh $117\mbox{k}$, set $-\alpha \Delta t^2 = 3.9$, $-\beta \Delta t = 1.9$ and $-\gamma = 0$, the computational time step is $\Delta t = T/720$ where $T = 1/f_e$ and the BDF1 time scheme is employed. 
Figure \ref{fig:Cd_osc_x} shows the time evolution of the drag coefficient. The amplitude of the signal is about 3.5 that results to be in very good agreement with that one computed in \cite{Shin2008,Dutsch1998}, that is 3.54.

\subsection{Eigenfrequencies of a rectangular cantilever beam coupled with a surrounding fluid} \label{sec:FSI}
In this subsection, we are going to test our algorithm for simulating FSI problems. The VBM, in its standard formulation as PI controller, was successfull used within FSI framework: see, e.g., \cite{Huang2007, Shin2010, Song2011, Qin2012, Uddin2015, Son2017}. Here, the main goal is to show the role played by the derivative controller term. We investigate the natural oscillation frequencies of a rectangular cantilever beam embedded in a fluid domain \cite{VanEysden2006, Sader1998}. The structural model equation used here, including inertia and bending effects, is

\begin{equation}\label{eq:structural1}
\rho_s A_t \dfrac{\partial^2 r_y}{\partial t^2} = -EI \dfrac{\partial^4 r_y}{\partial x^4} - A_t F_y, 
\end{equation}
where $x$ is the beam axis, $r_y$ is the transversal displacement, $\rho_s$ is the mass density, $A_t = bh$ is the transversal section with $b$ and $h$ the width and thickness, respectively, $E$ is the Young modulus, $I = bh^3/12$ is the moment of inertia, $F_y$ is the transversal component of the Lagragian forcing acted on the body by the surrounding fluid. The boundary conditions for the equation \eqref{eq:structural1} explored here refer to clamped-free configuration. They read

\begin{equation}\label{eq:structural2}
\begin{cases}
r_y\rvert_{x = 0} = {\dfrac{\partial r_y}{\partial x}}\bigg{\rvert}_{x = 0} = 0, \\ \\
\dfrac{\partial^2 r_y}{\partial x^2}\bigg{\rvert}_{x = L} = \dfrac{\partial^3 r_y}{\partial x^3}\bigg{\rvert}_{x = L} = 0.
\end{cases}
\end{equation}


We have chosen to implement the elasticity model \eqref{eq:structural1}-\eqref{eq:structural2} within the finite element C++ library \emph{deal.II}. We adopt linear finite elements for the space discretization. The second-order time derivative is discretized by a finite difference method based on a three-points centered scheme. The linear algebraic system associated is solved by using the Generalized Minimal Residual Method GMRM. The required accuracy is 1e-7 at each time step. The coupling process between flow solver and structural solver can be summarized as follows:

\begin{enumerate}
\item At $t^n$, we know $\u^{n}$ and $r_y^n$ and $r_y^{n-1}$. Then we calculate $F_y^n$ by eq. \eqref{eq:forc1};

\item We obtain $f_y^n$ by eq. \eqref{eq:forc5} and we solve the problem \eqref{eq:evolveFV-1.1}-\eqref{eq:Poisson} to obtain $\u^{n+1}$ and $p^{n+1}$;

\item We solve the problem \eqref{eq:structural1}-\eqref{eq:structural2} to obtain $r_y^{n+1}$.
\end{enumerate}

We compare our computational results with the following analytical predictions, here considered as \emph{true} solutions \cite{Sader1998}:

\begin{equation}\label{eq:analytical_frequency}
\dfrac{f_{fluid}}{f_{vacuum}} = \left(1 + \dfrac{\pi \rho b}{4 \rho_s h} \right)^{-1/2},
\end{equation}
where $f_{fluid}$ and $f_{vacuum}$ are the natural frequencies of the beam in fluid and vacuum, respectively. The natural frequencies in the vacuum, $f_{vacuum}$, are given by \cite{VanEysden2006, Sader1998}

\begin{equation}\label{eq:analytical_frequency_1}
f_{vacuum} = \dfrac{1}{2\pi}\dfrac{{K_n}^2}{L^2}\sqrt{\dfrac{EI}{\rho_s A_t}},
\end{equation}
where $K_n$ are the solutions of the following trascendental equation \cite{VanEysden2006, Sader1998}

\begin{equation}
1 + \cos(K_n L)\cosh(K_n L) = 0,
\end{equation}\label{eq:analytical_frequency_2}
where $L$ is the length of the beam. For further comparison, we will consider also the numerical data provided in \cite{Hengstler2013}, that performed two-way coupled 3D FSI simulations using the commercial simulation software ANSYS.

Geometrical and structural features of the beam are reported in Table \ref{tab:FSI1}. The fluid domain is a $1\times1$ rectangular box. The undisturbed vertical location of the beam is 0.5. The mesh consists of $10\mbox{k}$ cells. 20 cells are uniformly distributed along the beam and the remaining grid is stretched. The beam is discretized in 10 cells. The fluid is initially at the rest and do-nothing boundary conditions are imposed on all the boundaries. The beam is excited by a sinusoidal forcing applied at its free tip for 10 time steps with a total length of 0.2 ms \cite{Hengstler2013}. The final computational time is 0.2 s. We remark that in this case the choice of the gain coefficients can not be done based on the stability analysis reported in Sec. \ref{sec:stability_analysis} because the structural equation should be considered also. We set $-\alpha = 2e3$ and $-\beta = 7e2$ that have demonstrated to not lead to stability issues. Then, we carry out a sensitivity analysis with respect to $-\gamma$ ranging from 0 to 1.8. 


\begin{table}[h]
\centering
\begin{tabular}{ccccc}
\multicolumn{5}{c}{} \\
\cline{1-5}
$\rho$ [$Kg/m^3$] & $L [m]$ & $b [m]$ & $h [m]$ & $E [Pa]$ \\ 
\hline
2670 & 0.15 &  0.01 & 0.005 & 6.5e10 \\
\hline
\end{tabular}
\caption{Geometrical and structural features of the beam.}
\label{tab:FSI1}
\end{table}

\begin{table}[h]
\centering
\begin{tabular}{ccccc}
\multicolumn{5}{c}{} \\
\cline{1-5}
 & $K_n$ & Vacuum [Hz] \cite{VanEysden2006, Sader1998} & Air (\cite{VanEysden2006, Sader1998}/\cite{Hengstler2013}) [Hz] & Water (\cite{VanEysden2006, Sader1998}/\cite{Hengstler2013}) [Hz] \\ 
\hline
$f_1$ & 1.875  &  177 & 177/178 & 140/135 \\
$f_2$ & 4.694  &  1108 & 1108/1109 & 879/849 \\
\hline
\end{tabular}
\caption{First two bending modes computed by using the analytical predictions reported in \cite{VanEysden2006, Sader1998} (equation \eqref{eq:analytical_frequency}) and by numerical simulations reported in \cite{Hengstler2013} both for the air and the water.}
\label{tab:FSI2}
\end{table}

\begin{table}[h]
\centering
\begin{tabular}{cccccc}
\multicolumn{6}{c}{} \\
\cline{1-6}
Fluid & $-\gamma$ & $f_1$ (Hz) & $E_{f_1}$ (\%) & $f_2$ (Hz) & $E_{f_2}$ (\%) \\
\hline
Air & 0 & 175 & 1.1\% & 1080 & 2.5\%   \\
    & 0.5 & 175 & 1.1\% & 1080 & 2.5\%   \\
    & 1 & 175 & 1.1\% & 1080 & 2.5\%   \\
    & 1.8 & 175 & 1.1\% & 1080 & 2.5\%    \\ \hline
Water & 0 & 175 & 25\% & 1065 & 21.1\%    \\
    & 0.5 & 160 & 14.3\% & 985 & 12.1\%   \\
    & 1 & 150 & 6.7 \% & 925 & 5.2\%   \\
    & 1.8 & 140 & 0\% & 865 &  1.6\%  \\ 
\hline
\end{tabular}
\caption{First two bending modes computed by using the present algorithm both for the air and the water and percentage errors with respect to the analytical predictions reported in \cite{VanEysden2006, Sader1998} (equation \eqref{eq:analytical_frequency}).}
\label{tab:FSI3}
\end{table}

Table \ref{tab:FSI2} reports the first two bending modes, $f_1$ and $f_2$, of the cantilever beam both in air and in water computed by using the analytical formula \eqref{eq:analytical_frequency} \cite{VanEysden2006, Sader1998} and by the numerical simulations reported in \cite{Hengstler2013}. Table \ref{tab:FSI3} shows the first two bending modes obtained by using the FFT algorithm based on the average displacement computed by present simulations as well as associated percentage errors with respect to the analytical predictions \eqref{eq:analytical_frequency} \cite{VanEysden2006, Sader1998}. 
We observe that, for the air, the modes predicted with $-\gamma = 0$ are in very good agreement with the \emph{true} values and that the introduction of the derivative action does not affect the results. On the contrary, for the water, we note a large difference with respect to the corresponding \emph{true} values: in particular, the values obtained are practically the natural frequencies of the beam in the air. When $-\gamma \neq 0$, we observe a significant improvement of the results obtained that are closer and closer to the \emph{true} values at increasing of $-\gamma$. We speculate that, for the water, the derivative action is necessary for properly detecting the eigenfrequencies because the right amount of added mass should be taken into account. On the contrary, for the air, the amount of added mass is negligible and the method works well even without derivative action. To the best of our knowledge, this benchmark has been unexplored within the VBM framework. However, we observe that very similar problems such as flexible filaments \cite{Huang2007} and flapping dynamics of coupled flexible flags in a uniform flow \cite{Son2017} were successfull simulated with the standard formulation of the VBM. 

In conclusion, the present findings indicate that the introduction of the derivative action could play an important part in order to properly take into account the fluid-structure coupling.  


\section{Conclusions and Perspectives}\label{sec:conclusion}
We showed the effectiveness of a FV-based VBM solver in simulating flow problems involving fixed, moving and deforming bodies embedded in a fluid region. The interest in the Finite Volume approximation is due to the fact that, to the best of our knowledge, the application of the VBM in a FV framework has been unexplored. We proposed to modify the standard feedback forcing scheme of the VBM, based on a PI controller, with the introduction of a derivative action, in order to obtain a PID controller. We showed that the derivative action affects strongly the stability region and in different way, depending on the time scheme used. In order to showcase the features of our approach, we presented a computational study related to 2D flow past a moving/fixed rigid cylinder in several configurations and eigenfrequencies of a cantilever beam coupled with a surrounding fluid. We showed that when bodies at rest are considered, the accuracy and efficiency of the computation is essentially governed by the integral controller. On the other hand, when bodies in rigid motion are considered, proportional and derivative actions help to mitigate unphysical oscillations that affect the solution at large integral controller coefficients. Finally, in the FSI context, we showed the importance of the derivative action in order to obtain proper results when added mass is comparable to solid body mass. 

As a follow-up of the present work, we are going to carry out further investigations in order to better clarify the role played by the derivative action, especially within the FSI context. Moreover, we would like to use the present approach to simulate turbulent flows both in RANS and LES contexts. We are also interested in developing a VBM Reduced Order Model (ROM) within a Finite Volume framework.

\section*{Akwnoledgements}

We acknowledge the support provided by the European Research Council Executive Agency by the Consolidator Grant project AROMA-CFD “Advanced
Reduced Order Methods with Applications in Computational Fluid Dynamics” - GA 681447, H2020-ERC CoG 2015 AROMA-CFD and INdAM-GNCS projects.

\bibliographystyle{amsplain_mod}
\bibliography{bib/references.bib} 
\end{document}